\documentclass[10pt]{article}

\usepackage{amscd,amsmath, amssymb, fancyhdr, tikz, epsfig}

\usepackage[small,nohug]{fake_diagrams}
\diagramstyle[small]

\usepackage[backref=page]{hyperref}
\renewcommand*{\backrefalt}[4]{%
	\ifcase #1 (Not cited.)%
	\or        (Cited on page~#2.)%
	\else      (Cited on pages~#2.)%
	\fi}

\hypersetup{
	colorlinks   = true,
	citecolor    = magenta,
	linkcolor    = blue,
	urlcolor     = magenta	
}

\numberwithin{equation}{section}

\newcommand{\version}{version 6.0,\ \ June 5, 2026}

\def\eqref#1{(\ref{#1})}
\newcommand{\goth}{\mathfrak}

\newcommand{\arrow}{{\:\longrightarrow\:}}
\newcommand{\Z}{{\Bbb Z}}
\def\C{{\Bbb C}}
\newcommand{\R}{{\Bbb R}}
\newcommand{\Q}{{\Bbb Q}}

\newcommand{\6}{\partial}
\def\1{\sqrt{-1}\:}
\newcommand{\restrict}[1]{{\left|_{{\phantom{|}\!\!}_{#1}}\right.}}
\newcommand{\cntrct}                % contraction with a vector field
{\hspace{2pt}\raisebox{1pt}{\text{$\lrcorner$}}\hspace{2pt}}

\newcommand{\calo}{{\cal O}}

\renewcommand{\tilde}{\widetilde}
\renewcommand{\bar}{\overline}
\renewcommand{\phi}{\varphi}
\renewcommand{\epsilon}{\varepsilon}
\renewcommand{\geq}{\geqslant}
\renewcommand{\leq}{\leqslant}

\newcommand{\im}{\operatorname{im}}
\newcommand{\End}{\operatorname{End}}

\newcommand{\Id}{\operatorname{Id}}
\newcommand{\Gr}{\operatorname{Gr}}

\newcommand{\Vol}{\operatorname{Vol}}

\newcommand{\Hom}{\operatorname{Hom}}
\newcommand{\Sym}{\operatorname{Sym}}
\newcommand{\St}{\operatorname{St}}
\newcommand{\Discr}{\operatorname{Discr}}

\newcommand{\Aut}{\operatorname{Aut}}
\newcommand{\Alb}{\operatorname{Alb}}

\newcommand{\Lie}{\operatorname{Lie}}

\newcommand{\Jac}{\operatorname{Jac}}

\newcommand{\codim}{\operatorname{codim}}

\newcommand{\Tw}{\operatorname{Tw}}

\renewcommand{\Re}{\operatorname{Re}}
\renewcommand{\Im}{\operatorname{Im}}

\newcounter{Mycounter}[section]
\newcounter{lemma}[section]
\renewcommand{\thelemma}{{Lemma \thesection.\arabic{lemma}}}
\newcommand{\lemma}{%
    \setcounter{lemma}{\value{Mycounter}}
    \refstepcounter{lemma}
    \stepcounter{Mycounter}
    {\noindent \bf \thelemma:\ }}

\newcounter{claim}[section]
\renewcommand{\theclaim}{{Claim \thesection.\arabic{claim}}}
\newcommand{\claim}{%
    \setcounter{claim}{\value{Mycounter}}
    \refstepcounter{claim}
    \stepcounter{Mycounter}
    {\noindent \bf \theclaim:\ }}

\newcounter{sublemma}[section]
\newcounter{corollary}[section]
\renewcommand{\thecorollary}{{Corollary \thesection.\arabic{corollary}}}
\newcommand{\corollary}{%
    \setcounter{corollary}{\value{Mycounter}}
    \refstepcounter{corollary}
    \stepcounter{Mycounter}
    {\noindent \bf \thecorollary:\ }}

\newcounter{theorem}[section]
\renewcommand{\thetheorem}{{Theorem \thesection.\arabic{theorem}}}
\newcommand{\theorem}{%
    \setcounter{theorem}{\value{Mycounter}}
    \refstepcounter{theorem}
    \stepcounter{Mycounter}
    {\noindent \bf \thetheorem:\ }}

\newcounter{conjecture}[section]
\renewcommand{\theconjecture}{{Conjecture \thesection.\arabic{conjecture}}}
\newcommand{\conjecture}{%
    \setcounter{conjecture}{\value{Mycounter}}
    \refstepcounter{conjecture}
    \stepcounter{Mycounter}
    {\noindent \bf \theconjecture:\ }}

\newcounter{proposition}[section]
\renewcommand{\theproposition}
      {{Proposition \thesection.\arabic{proposition}}}
\newcommand{\proposition}{%
    \setcounter{proposition}{\value{Mycounter}}
    \refstepcounter{proposition}
    \stepcounter{Mycounter}
    {\noindent \bf \theproposition:\ }}

\newcounter{definition}[section]
\renewcommand{\thedefinition}
      {{Definition~\thesection.\arabic{definition}}}
\newcommand{\definition}{%
    \setcounter{definition}{\value{Mycounter}}
    \refstepcounter{definition}
    \stepcounter{Mycounter}
    {\noindent \bf \thedefinition:\ }}

\newcounter{example}[section]
\renewcommand{\theexample}{{Example \thesection.\arabic{example}}}
\newcommand{\example}{%
    \setcounter{example}{\value{Mycounter}}
    \refstepcounter{example}
    \stepcounter{Mycounter}
    {\noindent \bf \theexample:\ }}

\newcounter{remark}[section]
\renewcommand{\theremark}{{Remark \thesection.\arabic{remark}}}
\newcommand{\remark}{%
    \setcounter{remark}{\value{Mycounter}}
    \refstepcounter{remark}
    \stepcounter{Mycounter}
    {\noindent \bf \theremark:\ }}

\newcounter{problem}[section]
\newcounter{question}[section]
\newcommand{\pstep}{{\bf Proof. Step 1:\ }}
\newcommand{\proof}{{\bf Proof:\ }}

\makeatletter

\@addtoreset{equation}{section} \@addtoreset{footnote}{section}
\makeatother

\def\blacksquare{\hbox{\vrule width 5pt height 5pt depth 0pt}}
\def\endproof{\blacksquare}

\begin{document}

%%%%%%%%%%%%%%%%%%%%%%%%%%%%%%%%%%%%%%%%%%%%%%%%%%%%%%%%%%%%
\begin{center}
{\LARGE\bf
Sections of Lagrangian fibrations on 
holomorphic symplectic manifolds\\[4mm]
}
%%%%%%%%%%%%%%%%%%%%%%%%%%%%%%%%%%%%%%%%%%%%%%%%%%%%%%%%%%%%%

Fedor Bogomolov\footnote{Partially supported by a Simons Travel grant
and HSE University Basic Research Program}, 
Ljudmila Kamenova\footnote{Partially supported 
by a grant from the Simons Foundation/SFARI (522730, LK)}, 
Misha Verbitsky\footnote{Partially supported 
by FAPERJ SEI-260003/000410/2023 and CNPq - Process 310952/2
021-2. 

 }

\end{center}

{\small \hspace{0.10\linewidth}
\begin{minipage}[t]{0.85\linewidth}
{\bf Abstract.}  Let $M$ be a holomorphically
symplectic manifold, equipped with a Lagrangian fibration
$\pi:\; M \to X$. A degenerate twistor deformation 
(a more general version of it is called ``a Tate-Shafarevich twist'') is a family 
of holomorphically symplectic structures on $M$
parametrized by $H^{1,1}(X)$. All members of this
family are equipped with a holomorphic Lagrangian projection
to $X$, and their fibers are isomorphic to the fibers of $\pi$.
Assume that $M$ is a compact hyperk\"ahler manifold of maximal holonomy,
and the general fiber of the Lagrangian projection $\pi$ 
is primitive (that is, not divisible) in integer homology.
We also assume that $\pi$ has reduced fibers in codimension 1.
Then $M$ has a degenerate twistor deformation $M'$
such that the Lagrangian projection $\pi:\; M' \to X$
admits a meromorphic section. 
\end{minipage}
}

\tableofcontents

%%%%%%%%%%%%%%%%%%%%%%%%%%%%%%%%%%%%%%%%%%%%%%%%%%%%%%%%%%

\section{Introduction} 

%%%%%%%%%%%%%%%%%%%%%%%%%%%%%%%%%%%%%%%%%%%%%%%%%%%%%%%%%%

%%%%%%%%%%%%%%%%%%%%%%%%%%%%%%%%%%%%%%%%%%%%%%%%%%%%%%%%%%%%
\subsection{Holomorphic section of Lagrangian fibrations}
%%%%%%%%%%%%%%%%%%%%%%%%%%%%%%%%%%%%%%%%%%%%%%%%%%%%%%%%%%%%

A hyperk\"ahler manifold, for our present purposes,
is compact, K\"ahler, holomorphically symplectic manifold.
We say that $M$ is {\bf a hyperk\"ahler manifold of maximal holonomy, }
or {\bf IHS} (irreducible holomorphically symplectic) if 
$H^{2,0}(M)=\C$ and $\pi_1(M)=0$. The compactness and the maximal
holonomy condition is tacitly assumed throughout this paper. 

A {\bf Lagrangian fibration} on a hyperk\"ahler
manifold is a holomorphic map $\pi:\; M \arrow X$, 
with all fibers Lagrangian with respect to the
holomorphic symplectic form on $M$.
The manifold $X$ is biholomorphic to $\C P^n$
if it is smooth  (\cite{_Hwang:base_}).
With each Lagrangian fibration $\pi:\; M \arrow X$ we associate
a deformation $I_t$ of a complex structure on $M$,
parametrized by $H^2(X,\C)$ (\cite{_Verbitsky:degenerate_}).  
This deformation, called
{\bf the degenerate twistor deformation}, is properly
introduced in Subsection \ref{_Holo_Lagra_Subsection_}. It shares many 
properties with the twistor deformation, known
from hyperk\"ahler geometry; in particular,
its general fiber is non-algebraic.
Most importantly, each  member of this
family is equipped with a holomorphic Lagrangian projection
to $X$, and the fibers of this projection are isomorphic 
to the fibers of $\pi$.
The main result of the present paper is the
following theorem.

\hfill

%%%%%%%%%%%%%%%%%%%%%%%%%%%%%%%%%%%%%%%%%%%%%%%%%%%%%%%%%%%%
\theorem\label{_Sections_main_Theorem_}
Let $\pi:\; M \arrow X$ be a Lagrangian
fibration on a compact hyperk\"ahler manifold
of maximal holonomy, with $X=\C P^n$. 
Assume that the general fiber of the Lagrangian projection $\pi$ 
is primitive (that is, not divisible) in integer homology.
Assume, moreover, that the fibers of $\pi$ outside
of a codimension 2 subset in $X$ are reduced.
Denote by $(M,I_t)$ the degenerate twistor deformation,
with $t\in H^2(X, \C)$. Then for some $t_0\in H^2(X, \C)$,
the Lagrangian projection $\pi:\; (M,I_{t_0}) \arrow X$
admits a meromorphic section.

\proof \ref{_main_res_last_section_Theorem_}. \endproof

%%%%%%%%%%%%%%%%%%%%%%%%%%%%%%%%%%%%%%%%%%%%%%%%%%%%%%%%%%%%
\subsection{C-symplectic structures and
  degenerate twistor deformations}
%%%%%%%%%%%%%%%%%%%%%%%%%%%%%%%%%%%%%%%%%%%%%%%%%%%%%%%%%%%%

A superficially 
similar result was proven in \cite{_BDV:sections_}.
Recall that a closed complex-valued 2-form $\Omega$ on a smooth manifold
$M$, $\dim_\R M=4n$, is called {\bf C-symplectic} if 
$\Omega^n \wedge \bar \Omega^n$ is non-degenerate,
and $\Omega^{n+1}=0$. By 
\cite[Proposition 2.13]{_BDV:sections_},
for any C-symplectic form $\Omega$ on $M$ there exists
a unique complex structure operator $I$ such that
$\Omega$ is a holomorphic symplectic form on $(M,I)$.
The language of C-symplectic geometry is 
convenient in context of deformation theory, 
because it allows us to speak of deformations
of a holomorphically symplectic structure
without fixing the complex structure beforehand.
Using this language, one defines the degenerate twistor 
deformation as follows.

\hfill

\theorem\label{_dege_twistor_defo_Theorem_}
Let $M\stackrel{\pi}\arrow X$ be a Lagrangian fibration on a
holomorphically symplectic manifold $(X,\Omega)$, and
$\eta\in\Lambda^2(X)$ a closed
$(2,0)+(1,1)$-form on the base. Then the forms
$\Omega_t=\Omega+t\pi^*\eta$ on $X$ are C-symplectic, and
this deformation (called  {\bfseries\itshape the degenerate
  twistor deformation}) preserves the Lagrangian
fibration and the base.

\proof \cite[Theorem 3.3]{_BDV:sections_}. \endproof

\hfill

Using this observation and arguments of
linear-algebraic nature, Bogomolov, Deev and Verbitsky
proved the following theorem.

\hfill

\theorem\label{_Section_BDV_Theorem_}
Let $\pi:\; M \arrow X$ be a Lagrangian fibration 
on a holomorphic symplectic manifold $(M, \Omega)$,
and $S:\; X \arrow M$ a smooth ($C^\infty$) section.
Denote by $\eta\in \Lambda^2 (X)$ the form  
$-S^*(\Omega)$. Then $\Omega+ \pi^*\eta$ is
C-symplectic, and, moreover, the image 
of $S$ is holomorphic in the complex structure
$I_\eta$ induced by the C-symplectic form
$\Omega+ \pi^*\eta$.

\proof
\cite[Theorem 3.5]{_BDV:sections_}. \endproof

\hfill

\remark \ref{_Section_BDV_Theorem_}
 follows from \ref{_dege_twistor_defo_Theorem_}
and an observation of N. Hitchin 
(\cite[Proposition 1]{_Hitchin:Lagrangian_}), who proved 
that any submanifold $N\subset (M, \Omega)$
of a holomorphically symplectic manifold, 
which is Lagrangian with respect to $\Re\Omega$
and $\Im \Omega$, is holomorphic Lagrangian.
Indeed, the restriction of $\Omega+ \pi^*\eta$ to $S(X)$
vanishes by construction.

\hfill

After \ref{_Section_BDV_Theorem_}
was obtained, it seemed that existence of 
holomorphic sections of Lagrangian fibrations
is easily reachable by topological arguments.
After years of painstaking effort, this illusion shattered.
In fact, the topological arguments don't work
well even for K3 surfaces. In 
many years of trying we were unable to construct
a smooth section of an elliptic fibration on a 
K3 surface without using arguments from algebraic
geometry.
The aim of this paper is to develop this algebro-geometric
approach further, in order to solve the topological
problem of constructing a smooth section.

%%%%%%%%%%%%%%%%%%%%%%%%%%%%%%%%%%%%%%%%%%%%%%%%%%%%%%%%%%%%
\subsection{Shafarevich-Tate deformations and the relative Albanese scheme}
%%%%%%%%%%%%%%%%%%%%%%%%%%%%%%%%%%%%%%%%%%%%%%%%%%%%%%%%%%%%

The Shafarevich-Tate group has its origins in number theory,
but in the 1980-ies it appeared in the context of complex algebraic
geometry and topology, in the works of R. Friedman and 
J. Morgan on elliptic surfaces (\cite{_Friedman-Morgan-book_}).
Friedman and Morgan have defined the ``analytic Tate-Shafarevich group''
as the first cohomology of the sheaf of fiberwise automorphisms
of an elliptic surface over its base. Using an explicit computation
of this group, they proved that any elliptic surface can be
deformed to a surface with a section. 

The first attempt to generalize this notion to hyperk\"ahler
manifolds appeared in a work of E. Markman \cite{_Markman:Lagra_};
Markman produced a family of deformations of a hyperk\"ahler
manifold of K3${}^{[n]}$-type with a Lagrangian fibration, 
with the same fibers and the same base. Markman has explored
the K\"ahlerness of the deformation, and shown that
in some situations the K\"ahler property is preserved.
His construction was later generalized to all hyperk\"ahler
manifolds by A. Abasheva and V. Rogov, (\cite{_AbashevaRogov_,_Abasheva_II_}) who also
generalized his result about the K\"ahlerness.
Also, in \cite[Section 6]{_AbashevaRogov_}, the existence
of holomorphic sections of a degenerate twistor
deformation of a Lagrangian fibration was interpreted
as vanishing of a certain cohomology group related
to the Shafarevich-Tate group of this Larangian fibration.

The paper \cite{_Verbitsky:degenerate_} where the degenerate
twistor deformations were defined appeared at the same
time as Markman's \cite{_Markman:Lagra_}; however, it
was unknown what is the precise relation between these
two classes of deformations. The K\"ahlerness 
of the degenerate twistor deformations was finally proven in
\cite{_Soldatenkov_Verbitsky:degenerate_Kahler_}.

After the first version of this paper was published,
several important works dealing with the subject appeared.
In \cite{_Sacca:Lagrangian_}, G. Sacc\'a defined the relative Albanese variety,
a relative group scheme of fiberwise automorphisms of Lagrangian 
fibration on a hyperk\"ahler manifold without multiple fibers. 
The Shafarevich-Tate
deformations defined by Markman and Abasheva-Rogov naturally
appear as cohomology of the sheaf of local sections
of this relative group scheme. Another construction
of the relative Albanese variety appeared in \cite{_Kim:Neron_}
by Y.-J. Kim, who used a multi-dimensional version of the
N\'eron model.

One of the results of G. Sacc\`a \cite{_Sacca:Lagrangian_} is especially relevant to 
the present work. Using the minimal model program, she proves 
that the relative Albanese variety of a hyperk\"ahler manifold
$M$ admits a smooth hyperk\"ahler compactification. This gives
a manifold with a Lagrangian fibration admitting a holomorphic
section; however, it is not clear whether this new manifold
is deformationally equivalent (or even diffeomorphic) to $M$.

After the third version of this paper was published in arxiv.org,
J. Koll\'ar produced a version of our main result,\cite{_Kollar:sections_}.
In this paper, he shows that for any Abelian fiber space there exists
a Shafarevich-Tate twist which has a birational section. 
G. Sacc\`a  conjectured that any
Shafarevich-Tate twist of a hyperk\"ahler manifold $M$ is deformational
equivalent to $M$ (\cite[Remark 6.14]{_Sacca:Lagrangian_}), 
and proved that it has the same rational cohomology 
(\cite[Proposition 6.15]{_Sacca:Lagrangian_}).

%%%%%%%%%%%%%%%%%%%%%%%%%%%%%%%%%%%%%%%%%%%%%%%%%%%%%%%%%%%%
\subsection{Scheme of the proof}
%%%%%%%%%%%%%%%%%%%%%%%%%%%%%%%%%%%%%%%%%%%%%%%%%%%%%%%%%%%%

In this subsection, we give a rough outline of the arguments
used in this paper. 

An integer homology class, generally speaking, cannot be represented
by a smooth submanifold. However, for any smooth manifold $M$,
a class $z\in H_r(M, \Z)$ can be represented by 
a submanifold when $r \leq 5$. Such results 
were obtained for different dimensions and
classes of embeddings by  R. Thom, \cite{Thom}.
The strongest result
of this type is the following theorem.

\hfill

\theorem\label{_classes_realized_by_submanifold_}
For any orientable smooth $n$-manifold $V$, all elements
of the following integral homology groups can
be realized by orientable submanifolds:
$H_{n-1}(V,\Z)$, $H_{n-2}(V,\Z)$, $H_{i}(V,\Z)$
for all $i\leq 5$.

\proof \cite[Theorem II.27]{_NT:Cobordisms_1_}.
\endproof

\hfill

In conjunction with the theory of N\'eron models,
Thom's result is used to prove that, in the assumptions
of \ref{_Sections_main_Theorem_}, 
for any Lagrangian fibration $\pi:\; M \to X$,
there exists a smooth 2-dimensional submanifold $S\subset M$
such that the Lagrangian projection $\pi$ projects
$S$ to $\C P^1\subset X$ bijectively; in other words,
we obtain a smooth section of $\pi$ over a holomorphic
line $C=\C P^1 \subset M$ (\ref{_sections_from_CP^1_Theorem_}).

Taking a small neighbourhood of $C$, we can extend this
section to a small neighbourhood $U_C$, obtaining 
a map $\phi:\; U_C \arrow \pi^{-1}(U_C)$. Applying
\ref{_Section_BDV_Theorem_}, we immediately obtain 
a degenerate twistor deformation of the Lagrangian projection
$\pi^{-1}(U_C)\stackrel\pi\arrow U_C$ which makes this section holomorphic.

If we could extend this degenerate twistor deformation
to $M$, we would have been able to obtain a holomorphic section of
$\pi$ over a line $C \subset X$. However,
to extend the twistor deformation from $U_C$ to $X=\C P^n$, we need
to extend the closed $(2,0)+(1,1)$-form $\eta:= -\phi^*(\Omega)$
to a closed $(2,0)+(1,1)$-form on $\C P^n$. Such an extension,
generally speaking, has an obstruction, which is
an element in $H^1(\calo_{U_C})$ 
(Subsection \ref{_Dolbeault_classes_Subsection_}), 
and this group is infinite dimensional, as follows from 
a simple computation.

In Subsection \ref{_vanishing_Subsection_},
we compute this obstruction class, realizing it by a current of special type, which we call 
{\bf a Dolbeault current,}
and prove that it vanishes (\ref{_Dolbeault_vanishes_Theorem_}).

This allows us to construct a degenerate twistor
deformation $(M', \Omega')$ and a holomorphic section
of $\pi:\; M' \arrow X$ over a neighbourhood of a line $C\subset X=\C P^1$.
We notice that $C$ is an ample curve, that is, a curve with
an ample normal bundle. For an ample curve $C\subset X$, 
any holomorphic map $\phi:\; U_C\arrow Z$  from its neighbourhood $U_C$ 
to a complex variety $Z$ can 
be extended to a meromorphic map $\tilde \phi:\; X \arrow Z$ 
(\ref{_holo_for_Kahler_Theorem_}), and this finishes the proof.

%%%%%%%%%%%%%%%%%%%%%%%%%%%%%%%%%%%%%%%%%%%%%%%%%%%%%%%%%%%%%%%%%%%%%%%%
\subsection{Multiple fibers and primitivity}
\label{_multiple_fibers_intro_Subsection_}
%%%%%%%%%%%%%%%%%%%%%%%%%%%%%%%%%%%%%%%%%%%%%%%%%%%%%%%%%%%%%%%%%%%%%%%%

The assumption of primitivity of the general fiber of Lagrangian
projection is clearly necessary for existence of 
bimeromorphic section. In early versions of this
paper, we tried to treat it by using the notion of multiple
fibers in codimension 1: conjecturally, no Lagrangian
fibration on a hyperk\"ahler manifold may have multiple
fibers in codimension 1. Later, we became aware of the results 
of I. Hellmann \cite{_Hellmann:cone_},
who constructed a Lagrangian family with multiple fibers
(in codimension bigger than 1). After analyzing this example,
it became obvious that we use the primitivity of the homology
class of the general fiber of the Lagrangian projection. In the current
version of this paper, we made this depencency more explicit.

\hfill
 
We start by explaining the terminology, which  in some of published literature
might depend on conventions.

\hfill

\definition
Let  $\pi:\; M \arrow X$ be a proper holomorphic
map, $x\in X$ a point, $F_x:= \pi^{-1}(x)$ 
and $F_i$ its irreducible component, with 
(scheme-theoretical) multiplicity $\mu_i$.
Denote the greatest common divisor of $\mu_i$
by $\mu$. A fiber is {\bf multiple} if $\mu >1$.
A fiber $F_x$ is {\bf reduced} 
if $\mu_i=1$ for all $i$.
A fiber $F_x$ is {\bf has a  reduced component} 
if $\mu_i=1$ for at least one $i$.

\hfill

The following result should be well-known, but
we did not find it in the literature.

\hfill

%%%%%%%%%%%%%%%%%%%%%%%%%%%%%%%%%%%%%%%%%%%%%%%%%%%%%%%%%%%%
\lemma\label{_sections_iff_mult_1_Lemma_}
Let $f:\; X \arrow Y$ be a proper holomorphic
map of complex manifolds, and $y\in Y$ any point.
Then $f$ has a local section in a neighbourhood of $y$
if and only if the fiber $f^{-1}(y)$ has a component of
multiplicity 1. 

\hfill

\proof
Assume that $f^{-1}(y)$ has a reduced component $Z$, 
and consider a transversal to $Z$ local submanifold $U$
of complimentary dimension intersecting $Z$ in a smooth
point.  Since $U\cap f^{-1}(y)=1$,
for any $y'\in Y$ sufficiently close to $y$,
we also have $U\cap f^{-1}(y')=1$. The implicit
function theorem implies that the map
$y'\mapsto U\cap f^{-1}(y')$ is a holomorphic
map, defining a local section of $f$.

Conversely, assume that $f$ has a local section
$s:\; V\to X$ in a neighbourhood $V\ni y$. Let $s(V)=U$;
it is a smooth local submanifold of $X$.
The intersection $U\cap f^{-1}(y)$
has multiplicity 1, hence it belongs
to only one of the irreducible components $Z_1$ of $f^{-1}(y)$.
However, unless $Z_1$ is smooth in the intersection point
$U\cap Z_1$, the intersection $U\cap Z_1$ has multiplicity
$> 1$, which is impossible. Therefore, $Z_1$ is a
smooth submanifold of $X$ in a neighbourhood of the point
$U\cap Z_1$, and the multiplicity of this component is 1.
\endproof

\hfill

{\bf The discriminant} of a proper holomorphic map 
$\pi$ is the set of critical values of $\pi$.
By a theorem of Hwang and Oguiso (\cite{_Hwang_Oguiso:characteristic_}), 
the discriminant of a Lagrangian fibration is
a divisor in $X$. We say that the fibers of 
$\pi$ {\bf have a reduced component in codimension 1}
if the fiber $\pi^{-1}(s)$ is not multiple outside of a codimension
2 subset of $X$. It was conjectured that in codimension 1, 
the fibers of any Lagrangian fibration are actually reduced.
In an earlier version of this paper, we 
argued that the assumption of having a reduced component
for the fibers in codimension 1 is sufficient to have a
deformation admitting a bimeromorphic section.

\hfill

It is clear that when a bimeromorphic section exists,
the cohomology class of a general fiber should be primitive,
and this rules out the existence of multiple fibers. We still hope
that the bimeromorphic section exists when the fibers of 
$\pi$ have a reduced component in codimension 1.

Let $M \to \C P^n$ be a Lagrangian fibration on
a hyperk\"ahler manifold, and $C\subset \C P^n$ a general line.
We use the primitivity of a fiber in Section \ref{curve_sec},
when we argue for the existence of a 
2-dimensional submanifold $N$ in $M_C:=\pi^{-1}(C)$ 
such that $N\cap \pi^{-1}(pt)=1$, where $\cap$
means the intersection pairing in $H_*(M_C)$.
By R. Thom's theorem (\ref{_classes_realized_by_submanifold_}),
the existence of such $N$ is equivalent to 
primitivity of $\pi^*([H])$ in $H^2(M_C)$, where $[H]$
denotes the fundamental class of the hyperplane section.

Another place where the primitivity of the general fiber is used 
is Subsection \ref{_vanishing_Subsection_}. In the first
step of the proof of \ref{_Dolbeault_vanishes_Theorem_},
we invoke the same argument as used in Section \ref{curve_sec}
to produce the section over a curve. Here we apply it to
extend this section to a piecewise smooth multisection of
$\pi:\; M \to \C P^n$, in order to prove that the
Dolbeault current associated with the section of
$\pi$ in a neighbourhood of a line vanishes.

In \cite{_KV:primitive_}, Kamenova and Verbitsky
have  proven that $\pi^*([H])$ is primitive in $H^*(M)$,
assuming that $\pi:\; M \to \C P^n$ has no multiple fibers in  
codimension 1. However, it does not imply that $\pi^*([H])$ is primitive
in $H^*(M_C)$. Indeed, in \cite{_Hellmann:cone_}
I. Hellmann has shown that the Mukai system $M$ of rank two and genus two
on a K3 surface admits a Lagrangian fibration $\pi:\; M \to X$ with a
(very special) fiber of multiplicity 2. Therefore,
$\pi$ has no bimeromorphic sections. If $\pi^*([H])$
were primitive in $H^*(M_C)$, \ref{_Sections_main_Theorem_} 
would imply that $\pi$ has no multiple fibers.

However, existence of a section of $\pi:\; M \to \C P^n$
over a line is not sufficient. Indeed, in Hellmann's
example, such a section exists  (it exists for all
Beauville-Mukai Lagragian projections). For the benefit of
the reader, we explain this construction
now.\footnote{This 
argument is due to Giulia Sacc\'a.}

Let $M$ be a K3 surface.
The {\bf Beauville-Mukai manifold} is an irreducible
component of the moduli space of coherent sheaves
obtained as deformations of a line bundle of 
degree $d$ on a curve of genus $g$ in a K3 surface.

Fix a Beauville-Mukai system of all curves
of genus $g$ in $M$, equipped with a line bundle of degree $d$.
Appropriately compactified, this gives a hyperk\"ahler
manifold $X$ and a Lagrangian fibration $\pi:\; X \to Y$,
where $Y=\C P^n$ is the linear system of all curves of
genus $g$ on $M$, and $\pi$ is the forgetful map taking
the sheaf to its support variety with appropriate multiplicity.

Consider a rational pencil $C$ of
curves of genus $g$ passing through a fixed $z\in M$. This
is a rational curve on $Y$, and the Beauville-Mukai
Lagrangian family has a section $\sigma$ over $C$, which
takes a genus $g$ curve $S\subset M$ and picks the line bundle
$L={\cal O}(d[z])$; where $d[z]$ denotes a divisor
on $S$, obtained by taking $z$ $d$ times,

However, it is impossible to extend $\sigma:\; C \to X$ to 
a holomorphic section over a tubular neighbourhood of $C$.
A smooth section, which necessarily exists,
would violate \ref{_Dolbeault_vanishes_Theorem_}, that is,
the corresponding Dolbeault current is non-zero; in 
particular, it is impossible to extend the section
$\sigma$ to a rational section of $\pi:\; X \to Y$,
for any degenerate twistor deformation of $X$.

For Hellmann's example, the homology class
of the general fiber of the Lagrangian fibration is
manifestly non-primitive, hence a bimeromorphic
section cannot exist. We are grateful to C. Voisin, G. Sacc\'a and A. Abasheva
for this observation.

In \cite{_Kollar:sections_}, J. Koll\'ar proved a general
result which allows one to construct rational sections of 
abelian fibrations. His result is valid for any abelian fibration
with numerically trivial relative canonical bundle (this case
covers all hyperk\"ahler and Calabi-Yau manifolds). 
Instead of degenerate twistor transform (which makes no
sense outside of holomorphic symplectic geometry), 
Koll\'ar works with the Tate-Shafarevich twist of an
Abelian fibration $\pi:\; M \to X$; it is a fibration with the same base
and the same fibers, obtained from a 1-cocycle on $X$ with
coefficients in fiberwise automorphisms of $M$. 
Let $F$ be a general fiber of $\pi$ and $d:=\dim M-\dim X$.
Koll\'ar has shown that, given a homology class $\eta\in H_{2d}(M,\Z)$
such that $\eta \cap F=1$, there exists a Tate-Shafarevich twist 
$\pi_\eta:\; M_\eta\to X$ and a rational section of this fibration.

The existence of such $\eta$ implies non-existence of
multiple fibers. However, unlike the degenerate twistor
transform, the Tate-Shafarevich twist of $M$ does not need
to be deformationally equivalent (or even diffeomorphic) to $M$.

It is worth comparing our main result to results
of \cite{_Kollar:sections_}. In some of the earlier versions, we assumed
that the Lagrangian fibration has reduced fibers in 
codimension 1. In the current version, we drop this assumption,
following the idea of Koll\'ar (\ref{_sections_from_CP^1_no_Neron_Theorem_}).
However, our final statement is different from the one in \cite{_Kollar:sections_}, because we use the
degenerate twistor deformation, which does not change
the deformational type of a manifold, and  Koll\'ar
uses the Tate-Shafarevich twist, for which we do not know
whether it preserves the deformation class.
On the other hand, \cite{_Kollar:sections_}
applies to all fibrations with abelian fibers,
and our result works only for Lagrangian fibrations
on  hyperk\"ahler manifolds.

%%%%%%%%%%%%%%%%%%%%%%%%%%%%%%%%%%%%%%%%%%%%%%%%%%%%%%%%%%

\section{Preliminaries on holomorphic symplectic geometry} 
\label{_HK_Section_}

%%%%%%%%%%%%%%%%%%%%%%%%%%%%%%%%%%%%%%%%%%%%%%%%%%%%%%%%%%

%%%%%%%%%%%%%%%%%%%%%%%%%%%%%%%%%%%%%%%%%%%%%%%%%%%%%%%%%%%%%%%%%%%%%%%%
\subsection{C-symplectic structures}
%%%%%%%%%%%%%%%%%%%%%%%%%%%%%%%%%%%%%%%%%%%%%%%%%%%%%%%%%%%%%%%%%%%%%%%%

In this subsection, we introduce the formalism of
C-symplectic structures, used further on to explore
the degenerate twistor deformations. We follow 
\cite{_BDV:sections_} and \cite{_Soldatenkov_Verbitsky:Moser_}.

\hfill

\definition
Let $V$ be a real vector space. A {\bf complex structure operator} 
on $V$ is an element $I\in \Hom(V,V)$, satisfying $I^2=-\Id_V$. 

\hfill

\claim
The eigenvalues of $I$ are $\pm \1$. Moreover, $I$ is 
diagonalizable over $\C$. 
\endproof

\hfill

\definition
Let $V$ be a vector space, and $I\in \End(V)$ be a complex structure
operator. The eigenvalue decomposition $V\otimes_\R \C = V^{1,0}\oplus V^{0,1}$
is called the {\bf Hodge decomposition}; here $I\restrict{V^{1,0}}=\1\Id$,
and $I\restrict{V^{0,1}}=-\1\Id$.

\hfill

\remark
One can reconstruct $I$ from the space $V^{1,0}\subset V\otimes_\R \C$.
Indeed, take $V^{0,1}= \overline{V^{1,0}}$, and let $I$ act
on $V^{0,1}$ as $\1\Id$, and on $V^{0,1}$ as $-\1\Id$.
Since the operator $I\in \End(V\otimes_\R \C)$ commutes with  
complex conjugation, it is {\bf real}, that is, it preserves 
$V\subset V\otimes_\R \C$. 

\hfill

%%%%%%%%%%%%%%%%%%%%%%%%%%%%%%%%%%%%%%%%%%%%%%%%%%%%%%%%%%%%
\definition
An {\bf almost complex structure} on a real $2n$-manifold
$M$ is an operator $I\in \End(TM)$ satisfying $I^2=-\Id_{TM}$,
or equivalently, an $n$-dimensional subbundle 
$T^{1,0}M\subset TM\otimes_\R \C$ such that
$T^{1,0}M \cap \overline{T^{1,0}M}=0$. 
The almost complex structure is called {\bf integrable} 
(and $M$ is a {\bf complex manifold}) if $T^{1,0}M$ 
satisfies $[T^{1,0}M, T^{1,0}M]\subset T^{1,0}M$. 

\hfill

\remark Integrability implies existence of complex coordinates 
(by the classical theorem  of Newlander-Nirenberg).

\hfill

\definition
Let $(M,I)$ be a complex manifold, and let $\Omega\in\Lambda^2(M,\C)$
be a differential form. We say that $\Omega$ is {\bf non-degenerate} 
if $\ker \Omega\cap T_\R M =0$. We say that it is 
{\bf holomorphically symplectic} if it is non-degenerate, 
$d\Omega=0$, and $\Omega(IX,Y)= \1 \Omega(X,Y)$.

\hfill

\remark The equation $\Omega(IX,Y)= \1 \Omega(X,Y)$
means that $\Omega$ is complex linear 
with respect to the complex structure on $T_\R M$
induced by $I$. 

\hfill

\remark
Consider the Hodge decomposition $T_\C M= T^{1,0}M \oplus T^{0,1}M$
(the decomposition according to $I$'s eigenvalues). Since
$\Omega(IX,Y)= \1 \Omega(X,Y)$ and $I(Z)=-\1 Z$ for any
$Z\in T^{0,1}(M)$, we have  $\ker(\Omega)\supset T^{0,1}(M)$.
Since $\ker \Omega\cap T_\R M =0$, the real dimension of the kernel
is at most $\dim_\R M$, giving $\dim_\R \ker \Omega=\dim M$.
Therefore, $\ker(\Omega)= T^{0,1}M$.

\hfill

\corollary 
Let $\Omega$ be a holomorphically symplectic form on a
complex manifold $(M,I)$. 
Then $I$ is determined by $\Omega$ uniquely.

\hfill

\claim
Let $M$ be a smooth $2n$-dimensional manifold. Then
there is a bijective correspondence between the set
of almost complex structures and the set of sub-bundles
$T^{0,1}M\subset TM\otimes_\R\C$ satisfying $\dim_\C T^{0,1}M= n$
and $T^{0,1}M\cap TM=0$ (the last condition means that
there are no real vectors in $T^{1,0}M$, that is, that
$T^{0,1}M\cap T^{1,0}M=0$).

\hfill 

\proof Set $I\restrict{T^{1,0}M}=\1$ and
$I\restrict{T^{0,1}M}=-\1$.
\endproof

\hfill

\definition (\cite{_BDV:sections_})
Let $M$ be a smooth $4n$-dimensional manifold.
A closed complex-valued form $\Omega$ on $M$
is called {\bf C-symplectic} if $\Omega^{n+1}=0$
and $\Omega^{n}\wedge \bar\Omega^n$ is a non-degenerate
volume form. 

\hfill

\theorem
Let $\Omega\in \Lambda^2(M,\C)$ be a C-symplectic form,
and let $T^{0,1}_\Omega(M)$ be equal to $\ker \Omega$, where
\[
 \ker \Omega:= \{ v\in TM\otimes \C \ \ |\ \  i_v(\Omega)=0\},
\]
where $i_v$ denotes contraction with $v$.
Then $T^{0,1}_\Omega(M)\oplus 
\overline{T^{0,1}_\Omega(M)} = TM\otimes_\R \C$, hence
the sub-bundle $T^{0,1}_\Omega(M)$ defines an 
almost complex structure $I_\Omega$ on $M$. If, in addition, 
$\Omega$ is closed, then $I_\Omega$ is integrable, and
$\Omega$ is holomorphically symplectic on $(M, I_\Omega)$.

\hfill

\proof  The rank of $\Omega$ is $2n$ because $\Omega^{n+1}=0$
and $\Omega^n \neq 0$ everywhere. 
Then $\ker \Omega\oplus \overline{\ker\Omega}=T_\C M$.
The relation $[T^{0,1}_\Omega(M), T^{0,1}_\Omega(M)]\subset T^{0,1}_\Omega(M)$
follows from Cartan's formula 
\[
d\Omega(X_1,X_2, X_3)= \frac 1 6 
\sum_{\sigma\in \Sigma_3} (-1)^{\tilde \sigma}\Lie_{X_{\sigma_1}} 
\Omega(X_{\sigma_2}, X_{\sigma_3}) + 
(-1)^{\tilde \sigma}\Omega([X_{\sigma_1},X_{\sigma_2}], X_{\sigma_3})
\]
which gives, for all $X, Y \in T^{0,1}M$, and for any $Z\in TM$,
\[
d\Omega(X,Y,Z) = \Omega([X,Y], Z),
\]
implying that $[X,Y]\in T^{0,1}M$.
 \endproof

%%%%%%%%%%%%%%%%%%%%%%%%%%%%%%%%%%%%%%%%%%%%%%%%%%%%%%%%%%%%
\subsection{Holomorphic Lagrangian fibrations}
\label{_Holo_Lagra_Subsection_}
%%%%%%%%%%%%%%%%%%%%%%%%%%%%%%%%%%%%%%%%%%%%%%%%%%%%%%%%%%%%

In this subsection, we introduce holomorphic Lagrangian
fibrations, and collect their properties for further use.

\hfill

\definition
A {\bf holomorphic Lagrangian fibration} is a holomorphic map
$\pi:\; M \arrow X$ with $M$ holomorphically
symplectic and the smooth fibers are Lagrangian
with respect to the holomorphic symplectic form.

\hfill

\theorem (Arnold-Liouville) \label{AL}
Any smooth fiber of a proper Lagrangian fibration is a torus. 

\proof \cite{_Evans:Lagrangian_}. 
\endproof

\hfill 

\theorem (Matsushita, \cite{Mat}) 
Let $M$ be a hyperk\"ahler manifold of maximal holonomy, 
and let $\pi:\; M \arrow X$ be a surjective holomorphic map,
with $0 < \dim X < \dim M$. Then $\pi$ is a Lagrangian fibration. 

\hfill

\theorem (Hwang, \cite{_Hwang:base_}) 
Under the assumptions in Matsushita's theorem above, $X$ is 
biholomorphic to $\C P^n$ when it is smooth. 

\hfill

\conjecture Even under weaker assumptions, i.e., when $X$ is a normal 
variety, then $X$ is biholomorphic to $\C P^n$.

\hfill 

\remark Originally, Matsushita's and Hwang's theorems were only stated and 
proved in the projective setting. For a most general discussion and proof 
of these theorems, see \cite[Theorem 2.8]{KL} in the setting of singular
symplectic varieties, with a normal K\"ahler variety as the base. 

%%%%%%%%%%%%%%%%%%%%%%%%%%%%%%%%%%%%%%%%%%%%%%%%%%%%%%%%%%%%
\subsection{Degenerate twistor deformations}
%%%%%%%%%%%%%%%%%%%%%%%%%%%%%%%%%%%%%%%%%%%%%%%%%%%%%%%%%%%%

We round up the preliminaries by
exploring the basic properties of degenerate twistor deformations.

\hfill

\proposition  \label{_Hitchin_Lemma_Proposition_}
Let $Z\subset M$ be a submanifold of a holomorphically
symplectic manifold $(M, \Omega)$, $\dim M = 2\dim Z$. Assume
that $\Omega\restrict Z=0$. Then $Z$ is a complex submanifold.

\proof \cite[Proposition 1]{_Hitchin:Lagrangian_} or 
\cite[Proposition 2.7]{_BDV:sections_}. \endproof

\hfill

In \cite{_Verbitsky:degenerate_} properties of the degenerate
twistor deformation (\ref{_dege_twistor_defo_Theorem_}) 
on compact hyperk\"ahler manifolds 
were studied using the Hodge theory. 

\hfill

\remark
Let $\Tw_\eta(M) = M\times \C$ be equipped with a
complex structure $I_{tw}$ constructed as follows.
At any $(x, t)\in M \times \C$,
$I_{tw}$ acts on $T_{(x, t)}\Tw_\eta(M)= T_x M \oplus T_t\C$ 
as $I_{\Omega_{t\eta}}$ on $T_x M$ and as the standard complex
structure on $T_t\C$. It is not hard to see that $I_{tw}$
is integrable. Then $\Tw_\eta(M)$ is called the {\bf 
degenerate twistor space}. For a hyperk\"ahler manifold, 
it can be obtained as a 
limit of twistor spaces (\cite{_Verbitsky:degenerate_}), 
which explains the term. 

\hfill

\claim (\cite{_Verbitsky:degenerate_})\\
The map $\pi:\; M \arrow X$ is holomorphic with respect to $I_\eta$,
and the complex structures induced by $I$ and $I_\eta$ on the fibers of 
$\pi$ coincide. 
\endproof 

%%%%%%%%%%%%%%%%%%%%%%%%%%%%%%%%%%%%%%%%%%%%%%%%%%%%%%%%%%%%
\section{Canonical models for Lagrangian fibrations}
%%%%%%%%%%%%%%%%%%%%%%%%%%%%%%%%%%%%%%%%%%%%%%%%%%%%%%%%%%%%

%%%%%%%%%%%%%%%%%%%%%%%%%%%%%%%%%%%%%%%%%%%%%%%%%%%%%%%%%%%%
\subsection{Smooth torsor bundle associated with a
  Lagrangian fibration}
\label{_Canonical_torsor_Subsection_}
%%%%%%%%%%%%%%%%%%%%%%%%%%%%%%%%%%%%%%%%%%%%%%%%%%%%%%%%%%%%

We relate some results about N\'eron models
and their generalizations to Lagrangian fibered
hyperk\"ahler manifolds. We follow \cite{_KV:torsors_}.
The complex-geometric approach used in
\cite{_KV:torsors_} is based on the notion of
Kulikov model, discussed in 
\cite{_BLR:Nron_models_,_BHPT:abelian_CY_};
it was formally defined in \cite{_Kollar:Neron_}
as ``canonical model''.
 
The main result which we use appeared as a part of Theorem 1.6 
in \cite{_Kim:Neron_}. In \cite{_KV:torsors_} we gave an independent proof of the
following statement (see \cite[Theorem 1.6]{_Kim:Neron_}).

\hfill 

%%%%%%%%%%%%%%%%%%%%%%%%%%%%%%%%%%%%%%%%%%%%%%%%%%%%%%%%%%%%
\theorem\label{_main_torsor_Theorem_}
Let $M$ be a holomorphic symplectic manifold,
and $\pi:\; M \to X$ a Lagrangian fibration, that
is, a proper holomorphic map which has holomorphic
Lagrangian fibers and admits a relatively ample line 
bundle. Assume that $X$ is an open
subset of a projective manifold, and $M_i\subset M$
the set of all points where the differental
of $\pi$ has rank $\leq \dim X-i$.\footnote{In 
\cite{_KV:torsors_}, we show that $\codim \pi(M_i)\geq i$.}
Denote by
$X_0$ the set $X\backslash \pi(M_2)$,
and let $N\subset M$ be $\pi^{-1}(X_0) \backslash M_C$.
Then the natural projection $N \to X_0$
is equipped with a unique structure of 
a torsor of a smooth group scheme over the base $X_0$.

\proof
 \cite[Theorem 1.1]{_KV:torsors_}. \endproof

\hfill

Let us explain the terminology used in this statement.

\hfill

\definition
Let $f:\; X_0 \to Y$ be a smooth holomorphic submersion
(not necessarily proper), and $s:\; Y\to X_0$
a section. Consider a map $X_0\times_Y X_0\stackrel\mu \to X_0$
commuting with the projection to $Y$. We say that $(X_0,s, \mu)$
is {\bf a holomorphic group bundle} if the map $\mu$ is
submersive and defines a group structure on each fiber
$f^{-1}(y)$, for all $y\in Y$, and $s(y)$ is its identity element.

\hfill

By definition, a torsor over an abelian group
$G$ is a set where $G$ acts freely and transitively.
However, the torsor structure can be defined
without referring to an {\em a priori} given group structure as follows.

\hfill

\definition
{\bf An abelian torsor} is a set $S$ equipped
with a map $\mu:\; S\times S \times S \to S$ such that
\begin{description}
\item[(i)] for
any $e\in S$, the map $x \cdot y := \mu(x, y, e)$
defines a commutative and associative product on $S$
with $e$ an identity element 
\item[(ii)] 
 for any fixed $y \in S$,
the map $x\mapsto x\cdot y$ is bijective.
\end{description}

\remark
From the axioms of the abelian
torsor, it is clear that $\cdot$ defines a structure 
of an abelian group on $S$.

\hfill

\remark
The identity element axiom means that 
$\mu(x,y,y)=x$ for any $x$ and $y$. Given a group
structure on $S$, the torsor operation
$\mu$ can be defined  as $\mu(x, y,z)= x+y-z$.

\hfill

The most important property of abelian torsors is used
many times throughout this paper.

\hfill

%%%%%%%%%%%%%%%%%%%%%%%%%%%%%%%%%%%%%%%%%%%%%%%%%%%%%%%%%%%%
\claim\label{_sum_in_torsor_indep_Claim_}
Let $S$ be a set equipped with a structure of a torsor
of an abelian group $G$, and $a_1, ..., a_n, b_1, ..., b_{n+1}$
 points in $S$. Consider $e\in S$, and let $+$ denote the
abelian group structure on $S$ associated with the unit $e$.
Then the sum $\sum_{i=1}^{n+1}b_i -\sum_{i=1}^{n} a_i$
is independent from the choice of $e\in S$.

\hfill

\proof
Let us denote the group operation on $S$ by 
$x+_e y:=\mu(x, e, y)$,  and fix $u\in S$. 
Then $x+_u y= x+_e y-_e u$, and 
$$b_1+_u b_2+_u...+_u
b_{n+1}= b_1+_e b_2+_e...+_e b_{n+1}-_e  (n+1) u.$$
Then
\[ b_1+_u b_2+_u...+_u b_{n+1}-_u a_1-_u a_2-_u...-_u a_n= \]
\[b_1+_e b_2+_e...+_e b_{n+1}-_e a_1-_e a_2-_e...-_e a_n -_e  (n+1) u
+_e nu= \] 
\[b_1+_e b_2+_e...+_e b_{n+1}-_e a_1-_e a_2-_e...-_e
a_n
\]
(see also the proof of \ref{_sections_from_CP^1_Theorem_}, Step 2).
\endproof

%%%%%%%%%%%%%%%%%%%%%%%%%%%%%%%%%%%%%%%%%%%%%%%%%%%%%%%%%%%%
\subsection{N\'eron models for  Lagrangian fibrations}
\label{_Neron_Lagr_Subsection_}
%%%%%%%%%%%%%%%%%%%%%%%%%%%%%%%%%%%%%%%%%%%%%%%%%%%%%%%%%%%%

When a Lagrangian fibration $\pi:\; M \to X$
is equipped with a holomorphic section, the
torsor bundle defined in \ref{_main_torsor_Theorem_}
is actually a group bundle. In \cite{_Kim:Neron_},
this group bundle was interpreted as a part of the
multi-dimensional N\'eron model of $M$. Throughout
this paper, we often refer to this construction
as to the N\'eron model construction for the
Lagrangian fibration. Here is the relevant statement
from \cite{_KV:torsors_} formulated for Lagrangian
fibrations equipped with holomorphic sections.

\hfill

%%%%%%%%%%%%%%%%%%%%%%%%%%%%%%%%%%%%%%%%%%%%%%%%%%%%%%%%%%%%
\definition\label{_Discr_i_Definition_}
Let $M_i$ be the set of all points $m\in M$ such that
the differential $D\pi$ has rank $\leq \dim X-i$. 
The set $\pi(M_C)$ of singular values of a Lagrangian fibration map
$\pi:\; M \to X$ is called is {\bf discriminant},
denoted by $\Discr\subset X$. Let $\Discr_i:= \pi(M_{i+1})$.
Clearly, $\Discr=\Discr_0\supset \Discr_1\supset ...$.
By \cite[Claim 2.9]{_KV:torsors_}, $\codim_{\Discr} \Discr_i\geq i$.

\hfill

\remark
Note that $\Discr_i$ does {\em not} necesarily coincide with
the $i$-th singular stratum of $\Discr$. However, it has the
same dimension.

\hfill

The following theorem is a version of
\ref{_main_torsor_Theorem_}; in its earliest
form it appeared in \cite{_BHPT:abelian_CY_},
in the context of a ``Kulikov model''. 
Later, this construction was clarified in 
\cite{_Kollar:Neron_,_Sacca:Lagrangian_,_Kim:Neron_}.

\hfill

%%%%%%%%%%%%%%%%%%%%%%%%%%%%%%%%%%%%%%%%%%%%%%%%%%%%%%%%%%%%
\theorem\label{_canonical_model_group_Theorem_}
Let $M$ be a hyperk\"ahler manifold, 
$f:\; M \to X$ a Lagrangian fibration,
and $M_0\subset M$ its smooth locus.
Denote by $J(f)$ the manifold $M_0\cap f^{-1}(X\backslash \Discr_1)$.
Assume that $f$ admits 
a holomorphic section $s:\; X \to M_0$.
Then there exists a unique group bundle structure
on $f:\; J(f)\to X\backslash \Discr_1$ with the
following properties:
\begin{description}
\item[(i)] For each $y\in X\backslash \Discr_1$,
$s(y)$ is the identity element of the fiber $f^{-1}\cap J(f)$.
\item[(ii)] 
The natural action of the associated group bundle $J(f)$ 
on itself extends holomorphically to $f^{-1}(X\backslash \Discr_1)$.
\end{description}
\proof  \cite[Theorem 2.13]{_KV:torsors_}. \endproof

%%%%%%%%%%%%%%%%%%%%%%%%%%%%%%%%%%%%%%%%%%%%%%%%%%%%%%%%%%%%

\section{Constructing a smooth section over a curve} 
\label{curve_sec}

%%%%%%%%%%%%%%%%%%%%%%%%%%%%%%%%%%%%%%%%%%%%%%%%%%%%%%%%%%%%

%%%%%%%%%%%%%%%%%%%%%%%%%%%%%%%%%%%%%%%%%%%%%%%%%%%%%%%%%%%%
\subsection{Constructing a smooth section over a curve using N\'eron models} 
\label{_section_Nron_Subsection_}
%%%%%%%%%%%%%%%%%%%%%%%%%%%%%%%%%%%%%%%%%%%%%%%%%%%%%%%%%%%%

The main and only aim of this section is the following theorem.

\hfill

%%%%%%%%%%%%%%%%%%%%%%%%%%%%%%%%%%%%%%%%%%%%%%%%%%%%%%%%%%%%
\theorem\label{_sections_from_CP^1_Theorem_}
Let $\pi:\; M \to X$ be a hyperk\"ahler manifold equipped with a
Lagrangian fibration, with $X$ smooth (and hence,
by Hwang's theorem, isomorphic
to $\C P^n$), and $C\subset X$ a general line in $X=\C P^n$. 
Assume that the general fiber $F$ of the Lagrangian projection $\pi$ 
is primitive (that is, not divisible) in the integral homology. 
Assume, moreover, that the fibers of $\pi$ are reduced
away from a codimension $\geq 2$ subset of $X=\C P^n$.
Then
\begin{description}
\item[(i)] there exists a smooth section
$\sigma:\; C \arrow M$. 
\item[(ii)] Moreover,  this section
can be chosen homologous to a
homology class $\alpha \cap \pi^{-1}(C)\in H_2(\pi^{-1}(C), \Z)$,
where $\alpha\in H_{2n}(M, \Q)$ is a homology class such that
$\alpha \cap [F]=1$, where $[F]\in H_{2n}(M, \Z)$
denotes the homology class of a general fiber of $\pi$.
\end{description}

\remark
Since version 5 of this paper appeared in October 2025,
J\'anos Koll\'ar pointed out that the second claim of
this theorem is even more subtle than we originally thought. To fix
this gap, we had to add a subsection with the classification
of general singular fibers, due to Hwang and Oguiso
(Subsection \ref{_Hwang_Oguiso_Subsection_}). 
We resolve this new subtlety in Subsection
\ref{_extension_of_section_homology_Subsection_}. 

\hfill

\remark
The existence of local sections of $\pi:\; \pi^{-1}(C) \to C$ is needed
in \ref{_sections_from_CP^1_no_Neron_Theorem_}
(ii); it is also used in the original proof of 
\ref{_sections_from_CP^1_no_Neron_Theorem_} (i) 
based on the N\'eron models of the fibration
$\pi:\; \pi^{-1}(C) \to C$ (given below). In Subsection 
\ref{_curve_no_nron_Subsection_} we
give another version of 
\ref{_sections_from_CP^1_no_Neron_Theorem_} (i), 
which does not use N\'eron models and does not need this
assumption.

\hfill

We end this subsection with a proof of 
\ref{_sections_from_CP^1_Theorem_} (i).
In Subsection \ref{_extension_of_section_homology_Subsection_}, 
we explain why  \ref{_sections_from_CP^1_Theorem_} (ii)
is non-trivial and prove it, and
in Subsection \ref{_curve_no_nron_Subsection_}, we give another
proof of \ref{_sections_from_CP^1_Theorem_} (i).

\hfill

{\bf Proof of \ref{_sections_from_CP^1_Theorem_} (i), step 1:}
By Bertini's theorem, the preimage $M_C:=\pi^{-1}(C)$ of $C$ is smooth.
Denote by $F_1\subset M_C$ the preimage of a general point.
We are going to show that its fundamental
class $[F_1]\in H_{2n}(M_C)$ is primitive in the integral cohomology,
that is, is not divisible by a non-invertible integer.
Let $F:= \pi^{-1}(pt)$ be a general fiber of $\pi:\; M \arrow \C P^n$,
and $[F]\in H_{2n}(M,\Z)$ its homology class. Since $[F]$ is
primitive, there exists a homology class $\alpha\in H_{2n}(M,\Z)$ 
such that $[F]\cap \alpha=1$. Consider the class 
$\beta:=\alpha\cap [M_C]\in H_2(M_C, \Z)$; it can be realized as
an intersection of $\alpha$ and $M_C$, if the chains defining
$\alpha$, are transversal to $M_C$. By construction, this 
class intersects $[F]$ in one point, if counted with multiplicities,
which implies that $[F]\cap \beta=1$. 

Essentially, this argument boils down to the following: if $\pi$ has a 
singular (as in ``singular chain'') multisection $\alpha$ of total
multiplicity 1, then the intersection of $\alpha$
with $M_C$ is a singular chain which gives a multisection
of the map $\pi:\; M_C \arrow C$ of total multiplicity 1.

\hfill

{\bf Step 2:}
By the classical result of R. Thom 
(\ref{_classes_realized_by_submanifold_}), $\eta_1\in  H_2(M_C, \Z)$
can be represented by an embedded 2-submanifold $S\subset M_C$.
We can assume that $S$ does not meet the singular locus of $\pi:\; M_C\arrow C$
by dimensional reasons. Then, using the transversality
theorem (\cite{_EM:h-principle_}), 
we can also assume that $S$  intersects the fibers of $\pi$ in finitely many
points. Since $\eta_1 \cap F=1$, where $F$ is a fiber,
these points add up to 1, if counted with multiplicities. 
Let $T$ be a smooth fiber of $\pi:\; M_C \arrow C$;
by Arnold-Liouville's theorem (\ref{AL}),  $T$ is a torus. 
Denote the points of $S\cap T$ by $x_1, ..., x_{k+1}, y_1, ..., y_k$,
where $x_i$ are points of positive multiplicity, and $y_i$
of negative multiplicity. Fix a reference point $e\in T$.
This turns $T$ into a compact abelian group $T_0$.
Consider the point
\begin{equation}\label{_sum_with_e_depencence_Equation_}
z:=e+ \sum_{i=1}^{k+1}(x_i-e) - \sum_{i=1}^{k}(y_i-e)
\end{equation}
By \ref{_sum_in_torsor_indep_Claim_}, 
$z$ is independent of the choice of $e\in T$.

We refer to $a$ as to ``the average'' of 
$\{x_i, -y_i\}$, and we call this approach to constructing
a section ``the averaging''. \endproof

\hfill

\remark
This result was conjectured in \cite[Remark 4.2]{_KV:primitive_}.

%%%%%%%%%%%%%%%%%%%%%%%%%%%%%%%%%%%%%%%%%%%%%%%%%%%%%%%%%%%%
\subsection{Homology class of the averaging of a multisection} 
\label{_averaging_homology_Subsection_}
%%%%%%%%%%%%%%%%%%%%%%%%%%%%%%%%%%%%%%%%%%%%%%%%%%%%%%%%%%%%

In the proof of \ref{_sections_from_CP^1_Theorem_} (i), Step 2 
(Subsection \ref{_section_Nron_Subsection_}), we  used
the averaging procedure to produce a smooth section of an
abelian fibration $\pi:\; M_C\to C$ over a curve. Since the first version
of this paper was published, the same procedure was
used as the key argument in \cite{_KV:torsors_,_CKV:multisections_}.
The averaging takes a second cohomology class $[S]$ 
complementary\footnote{Here ``complementary'' means that
$[F]\cap [S]=1$.} to 
the homology class of a general fiber $F$ of $\pi$, represents
it by an embedded surface $S$, and averages
(indeed, adds together, counting multiplicities)
the intersection points of $S$ with each fiber of $\pi$,
obtaining a smooth section $S_0$. This
argument uses the fiberwise torsor structure on each fiber,
established in \cite{_Kollar:Neron_} using N\'eron models
(see \cite{_KV:torsors_} for more detail).

The result of this procedure depends on the choice of
the 2-submanifold $S$ representing the homology class $[S]$,
and $S_0$ does not need to be homologous to $S$.
In the present subsection we express the homology 
class of $S_0$ in terms of the homology class of $S$
and the characteristic leaves of the special fibers
(we define the characteristic leaves in Appendix, Section
\ref{_HO_Appendix_Section_},
following Hwang and Oguiso).

Let $V_1, V_2$ be two irreducible components of the same special
fiber $F_s$, and $W$ a connected component of their intersection.
Define {\bf the intersection sleeve} associated with $V_1, V_2, W$
as the homology class of an embedded 2-dimensional
sphere $S_{V_1, V_2, W}$, constructed as follows.

Let $a_1, a_2$ be two points in a smooth part
of $V_1, V_2$, situated near $W$, and $\Delta\subset C$
a small disc centered in $s=\pi(F_s)$.
Then $\pi^{-1}(\6\Delta)$ is a smooth submanifold
obtained as a torus bundle over a circle $\6 \Delta$,
and $\zeta$ a path connecting $a_1$ and $a_2$ in $F_s$.
Denote by $u$ the retraction map from $\pi^{-1}(\Delta)$ to
$F_s$, and let $S_{V_1, V_2, W}$ be the union of the cylinder
$u^{-1}(\zeta) \cap \pi^{-1}(\6\Delta)$
and two discs at both ends, $u^{-1}(a), u^{-1}(b)$.
\begin{center}
\epsfig{file=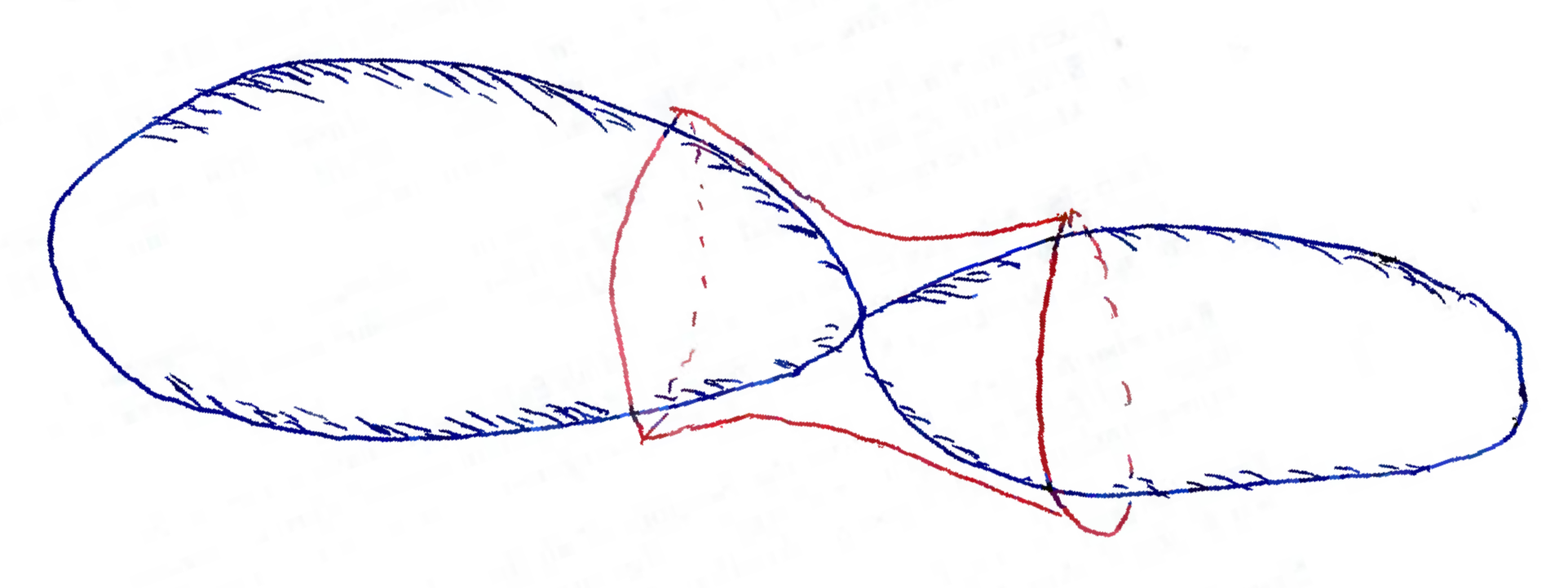,width=0.55\linewidth}\\
{\scriptsize The intersection sleeve (in red)}
\end{center}

Given a section $\sigma:\; C \to M_C$
of $\pi:\; M_C\to C$, intersecting
a special fiber $F_s$ in a component $V_1$,
we can always make a surgery by replacing the
disc $\sigma(\Delta)=u^{-1}(a)$ by the union of the cylindrical part
of $S_{V_1, V_2, W}$ and the disk $u^{-1}(b)$.
The new section is homologous to the
sum of the old one and the homology 
class of $S_{V_1, V_2, W}$. When $V_i$ has multiplicity $\mu_i$,
the disks $u^{-1}(a)$, $u^{-1}(b)$ intersect $V_1$, $V_2$
in $\mu_1$ and $\mu_2$ number of points respectively.

\hfill

%%%%%%%%%%%%%%%%%%%%%%%%%%%%%%%%%%%%%%%%%%%%%%%%%%%%%%%%%%%%
\lemma\label{_H_2_M_C_Lemma_}
Let $\pi:\; M \to \C P^n$ be a Lagrangian fibration,
$\pi:\; M \to C$ its restriction to a general
curve $C \subset \C P^n$, and $S_{V_1, V_2, W}$
the intersection sleeve constructed earlier in
this section. Then $S_{V_1, V_2, W}$ is homologous
to $[L_1]-[L_2]$, where $L_i$ are homology
classes of characteristic lines in $V_i$.

\hfill

\pstep
We use the notation introduced earlier when we constructed
$S_{V_1, V_2, W}$.  Let $\Delta_C$ be a neighbourhood of 
$s\in C$. The group $H^2(\pi^{-1}(\Delta_C))$
was computed in \ref{_second_cohomo_fiber_Proposition_} in
the Appendix.
We are going to pair $S_{V_1, V_2, W}$
with the generators $\alpha_i \in H^2(\pi^{-1}(\Delta_C))$
and conclude that 
$\langle \alpha_i, S_{V_1, V_2, W}\rangle= 
\langle \alpha_i, \mu_1[L_1]-\mu_2[L_2]\rangle$.
The space $H^2(\pi^{-1}(\Delta_C))= H^2(F_s)$
was described explicitly in
\ref{_second_cohomo_fiber_Proposition_}.

Consider the Albanese map
$\Alb:\; F_s \to V$ from the special 
fiber of $\pi$ to the torus $V$
(we give a description of this map in
the proof of \ref{_second_cohomo_fiber_Proposition_}, 
Step 2 in the Appendix).

Clearly, the image $\Alb(S_{V_1, V_2, W})$ is
contractible. Therefore, $S_{V_1, V_2, W}$ paired
with $\Alb^*(H^2(V))$ vanishes. The same
argument also implies that the pairing
of $S_{V_1, V_2, W}$ with the component
$H^1(V)$ of \eqref{_second_cohomo_when_circle_Equation_} 
vanishes.   It remains to compute the
pairing of $S_{V_1, V_2, W}$
with the classes $l_i\in H^2(F_s)$
constructed in \ref{_second_cohomo_fiber_Proposition_}.

\hfill

{\bf Step 2:}
We will use Poincar\'e duality
on a non-compact smooth manifold
obtained as a neighbourhood of the special fiber $F_s$
in $M_C:= \pi^{-1}(C)$, where $C\subset \C P^n$
is a general line.
Let $V_1,..., V_m$ be the irreducible
components of $F_s$, and $L_i \subset V_i$
the corresponding characteristic leaves.
The homology classes of $L_i$ and $V_i$ have
complementary dimension, because $\dim M_c =n+1$,
and the intersection  matrix $(L_i\cap V_j)$ is identified with the
Cartan matrix of the corresponding affine
Dynkin diagram. This Cartan matrix is 
degenerate, with radical of rank 1
obtained as $\sum m_i[L_i]$,
where $m_i\in \Z$ are coefficients
written on the nodes in the graphs
in Subsection \ref{_Hwang_Oguiso_Subsection_}.
These numbers are equal to the multiplicities
of the corresponding irreducible components.

A characteristic leaf meets each
irreducible component of $F_s$
if and only if the corresponding nodes are
connected, and the intersection index is
equal to the multiplicity. Therefore,
whenever $i\neq j$, $V_j \cap L_i= m_j$ if the
line $L_i$ intersects with $V_j$, and zero otherwise. 
To compute $V_i\cap L_i$, we notice that
$L_i\cap F_s=0$, because $F_s$ is homologous
to all other fibers of $\pi:\; M_C \to C$.
Therefore, $L_i\cap \sum_{j=1}^m m_j V_j=0$,
which implies $L_i\cap V_i= - m_{d_1}- ...- m_{d_k}$,
where $d_j$ are the numbers of the nodes connecting to
$L_i$, and $m_{d_j}$ their multiplicities.

If we prove that $S_{V_1, V_2, W}\cap V_j=([L_1]-[L_2])\cap V_j$
for any $V_j$, it will follow that $S_{V_1, V_2, W}-([L_1]-[L_2])$
belongs to the kernel of the pairing between the space 
$\langle L_i\rangle$ generated by the homology classes $L_i$ and 
$\langle V_j\rangle$, and this kernel is one-dimensional
and generated by the sum $\sum m_i L_i$.
To check that the difference $S_{V_1, V_2, W}-([L_1]-[L_2])$ vanishes,
we notice that $\sum m_i L_i$ is effective, hence it intersects
with the K\"ahler form $\omega$ positively, but 
$\langle S_{V_1, V_2, W}, \omega\rangle=0$. To prove
\ref{_H_2_M_C_Lemma_}, it remains to show that
$S_{V_1, V_2, W}\cap V_j=([L_1]-[L_2])\cap V_j$ for all $V_j$.

\hfill

{\bf Step 3:} Let $V_j':= m_j^{-1} V_j$.
Then $V_1' \cap S_{V_1, V_2, W}=1$,
$V_2' \cap S_{V_1, V_2, W}=-1$, 
$L_1 \cap V'_2= L_2 \cap V'_1=1$,
and $V_i' \cap L_i=-1 - \sum_j \frac{m_{d_j}}{m_i}$ 
where the sum is taken over all nodes connected to $V_i$.
Looking at the diagrams in Subsection \ref{_Hwang_Oguiso_Subsection_},
we see that this number is always -2. This is not
very surprising, because the affine Dynkin diagrams
are intersection forms for irreducible
components of the special fibers in elliptic
K3 surface. Unless the special fiber has only
one irreducible component (in which case
it is a nodal rational curve), each of the components of the
special fiber is a smooth rational curve which
has self-intersection -2.

Now, the intersection form between $V_i'$ and $L_j'$
is represented by the matrix $\begin{pmatrix} -2 &1\\1& -2\end{pmatrix}$,
and $V_1' \cap S_{V_1, V_2, W}=1$,
$V_2' \cap S_{V_1, V_2, W}=-1$ implies that 
the homology class of $S_{V_1, V_2, W}$ is in the orthogonal
complement of $V_1'+ V_2'$, that is, proportional to
$L_1-L_2$. Since $S_{V_1, V_2, W}\cap V_1= m_1$,
this also implies that $S_{V_1, V_2, W}$ is 
homologous to $L_1-L_2$.
\endproof

\hfill

The main result of this subsection is the following
theorem.

\hfill

%%%%%%%%%%%%%%%%%%%%%%%%%%%%%%%%%%%%%%%%%%%%%%%%%%%%%%%%%%%%
\theorem \label{_averaging_class_Theorem_}
Let $\pi:\; M \to \C P^n$ be a Lagrangian fibration
on a hyperk\"ahler manifold, and $C\subset \C P^n$
a general curve. Consider the fibration
$\pi:\; M_C \to C$, where $M_C:= \pi^{-1}(C)$,
let $S\subset M_C$ be a 2-dimensional submanifold
such that its homology intersection with
the general fiber of $\pi$ is 1, 
and let $S'$ be the section obtained from
$S$ by the averaging procedure 
(\ref{_sections_from_CP^1_Theorem_}, Step 2). Assume that the homology
intersection of $S$ with $V_i\in \pi^{-1}(S)$ 
is $a_i \in \Z$, and $[S']\cap V_i=b_i$;
by construction, $b_i=0$ for all $i$
except one. Renumbering the components,
we can assume that $b_1=1$.
Then the homology class $[S']$ is equal to 
$[S] - \sum a_i [L_i] + [L_1]$.

\hfill

\pstep
Using  \ref{_canonical_model_group_Theorem_},
we remove a codimension 2 subvariety from $M_C$,
obtaining a group fibration over $C$, with
$S'$ as its unit section.
Let $z\in C\cap D$, and denote by $V_1$ 
the unit component of $\pi^{-1}(z)$.
Using the same surgery as used in the definition in
the intersection sleeve, we can always modify $S$
by adding a collection of intersection sleeves,
and obtain a new multisection, denoted $S_1$,
which intersects $\pi^{-1}(z)$
in $V_1$ only. By construction, the averaging of
$S_1$ is homologous to $S'$; indeed, the averaging
can be obtained by a continuous deformation $S_t$,
as long as $S_t$ does not intersect the points
in $V_i\cap V_j$, where $V_i$ are irreducible
components of special fibers. On the other hand,
$[S'] = [S_1]= [S]+ \sum_{i_j} S_{V_{i_j}, V_{i_j'}, W_{i_j}}$,
where $\{S_{V_{i_j}, V_{i_j'}, W_{i_j}}\}$
denotes the set of all intersection sleeves
used in the surgery. However, the homology
class of each sleeve is $[L_{i_j}]-[L_{i_j'}]$,
hence $[S'] = [S]+ \sum_i c_i [L_i]$,
for some integer numbers $c_i$ depending on 
the manner of the surgery.

\hfill

{\bf Step 2:}
We know the intersections of $S$ and $S'$
with the irreducible components $V_i$
of the fibers, and we know that $[S'] = [S]+ \sum_i c_i [L_i]$.
To finish the proof of \ref{_averaging_class_Theorem_},
it remains to show that these two bits of information
are sufficient to recover the difference $[S]-[S']$.
However, the intersection pairing between
the homology group generated by $[L_i]$
and the homology group generated by $[V_i]$ has
a kernel, as explained below (see also the proof of
\ref{_H_2_M_C_Lemma_}), and we need to 
account for this as well.

Let $W\subset H_2(M_C, \R)$ be the vector space
generated by $[L_i]$, and $b$
be the intersection pairing 
between $W$ and the subspace $W_1$ in $H_{2n}(M_C, \R)$
generated by the irreducible components of $\pi^{-1}(z)$.
Consider the map $W\stackrel \rho \to H_{2n}(M, \R)$
taking each $L_i$ to the corresponding component $V_i$ of the fiber
$\pi^{-1}(z)$. Then $b([L_i], [V_j])= [V_i]\cap[V_j]$,
where $\cap$ denotes 
the intersection pairing between these
components in $H^{2n}(M, \R)$.
From the description obtained in \ref{_HO_classification_Theorem_},
it is apparent that this pairing is identified with 
the Cartan matrix of the corresponding affine Dynkin diagram.
Its radical $r$ is generated by the sum of all $[V_i]$ with 
coefficients equal to the multiplicities.
Denote by $W'\subset W$ the subspace
generated by $[L_i]-[L_j]$.

In Step 3 we are going to show that $W'\not\ni r$.
Since $[S]-[S']\in W'$, this implies that
the element $[S]-[S']$ is uniquely
determined by its intersection pairing
with all $V_i$. Clearly, 
$[S']\cap [V_j]=([S] - \sum a_i [L_i] + [L_1])\cap [V_j]$
because both sides of this equation
are equal to $\delta_{1j}$.\footnote{By definition,
$\delta_{1j}$ is $1$ when $j=1$ and 0 when $j\neq 1$.}
This proves \ref{_averaging_class_Theorem_}
modulo $W'\not\ni r$.

\hfill

{\bf Step 3:} It remains to prove that $W'\not\ni r$,
in the assumptions and notations of Step 2. 
Consider the functional $W\to \R$
taking $[L_i]$ to its multiplicity $m_i$.
This functional vanishes on $W'$ and does not
vanish on $r$ because the sum of all multiplicities
of $r= \sum m_i [L_i]$ is positive.
\endproof

\hfill

The same argument also proves the following corollary.

\hfill

%%%%%%%%%%%%%%%%%%%%%%%%%%%%%%%%%%%%%%%%%%%%%%%%%%%%%%%%%%%%
\corollary\label{_modulo_sum_determined_by_prod_Corollary_}
In the assumptions of \ref{_averaging_class_Theorem_},
let $l\in H_2(M_C, \Z)$ be a homology class
obtained as a linear combination of the 
homology classes of characteristic lines
in $\pi^{-1}(z)$, where $z\in C\cap D$.
Then $l$  is uniquely, modulo $\sum m_j l_j$, determined by
the intersection numbers $l\cap V_i$,
where $V_1, ..., V_s$ are irreducible
components of $\pi^{-1}(z)$.
\endproof

\hfill

\ref{_averaging_class_Theorem_} is interesting in itself,
but we are not going to use an explicit expression of the
class $[S]-[S']$ obtained there. In the rest of this paper,
we shall only use the following corollary
of \ref{_averaging_class_Theorem_}.

\hfill

%%%%%%%%%%%%%%%%%%%%%%%%%%%%%%%%%%%%%%%%%%%%%%%%%%%%%%%%%%%%
\corollary\label{_corollary_S'_independent_Corollary_}
Let $\pi:\; M \to \C P^n$ be a Lagrangian fibration
on a hyperk\"ahler manifold, $C\subset \C P^n$
a general curve, and $M_C:=\pi^{-1}(C)$.
Consider a homology class $[S]\in H_2(M_C, \Z)$
such that $[S]\cap [F]=1$, where $[F]$ is the
homology class of the general fiber of $\pi:\; M_C\to \C P^1$.
Let $S'$ be the smooth section of $\pi:\; M_C \to C$
obtained from $S$ by the averaging procedure (\ref{_sections_from_CP^1_Theorem_},
Step 2). Then the homology class of $S'$ in $H_2(M_C, \Z)$
is uniquely determined by $[S]$.
\endproof

\subsection{Homology classes of sections of abelian fibrations
over a line extended to the ambient manifold} 
\label{_extension_of_section_homology_Subsection_}
%%%%%%%%%%%%%%%%%%%%%%%%%%%%%%%%%%%%%%%%%%%%%%%%%%%%%%%%%%%%

In \ref{_sections_from_CP^1_Theorem_} (i)
we proved that there exists a smooth section
of a Lagrangian fibration restricted to a line 
$\C P^1 \subset \C P^n$. This section is obtained
by taking a multisection of $M_C\arrow \C P^1$
and applying the averaging procedure as explaned
in Subsection \ref{_section_Nron_Subsection_}, Step 2.

After an earlier version of this paper was finished,
J. Koll\'ar contacted us remarking that this
procedure might change the homology class of the
section (\ref{_averaging_changes_class_Example_}). 
This concern is well founded, and we address it in the
present subsection.

In the paper \cite{_KS:Sections_},
J. Koll\'ar and G. Sacc\'a 
studied universal compactified Jacobian $\pi:\; M \to \C P^n$
of a family of curves on a projective surface.
This universal compactified Jacobian is a projective
variety fibered over a line system of divisors
in $S$ (which are curves) with a general fiber
over $C\subset S$ identified with the Jacobian $\Jac(C)$.
A special case of such an abelian fibration is the
classical Beauville-Mukai system on a K3 surface
(Subsection \ref{_multiple_fibers_intro_Subsection_}).
Fix a general line $C=\C P^1\subset \C P^n$,
and let $M_C:= \pi^{-1}(C)$. J. Koll\'ar and G. Sacc\'a
study the restriction map  $r:\; H^{2g}(M, \Z)\to H^{2g}(M_C, \Z)$,
and show that is not surjective. Indeed, for any
holomorphic section $s$ of $M_C \to C$, such that its 
fundamental class belongs to $\im (r)$, $s$ is necessarily 
a zero section of the fibration $M_C \to C$, understood
as a group scheme over $C$.

In an email letter, J. Koll\'ar gave the following 
remarkable example.

\hfill

%%%%%%%%%%%%%%%%%%%%%%%%%%%%%%%%%%%%%%%%%%%%%%%%%%%%%%%%%%%%
\example\label{_averaging_changes_class_Example_}
Let $\pi:\; X \to \C P^1$ be an elliptic fibration on a K3 surface.
Assume that a special fiber of $\pi$ is a collection
of at least 3 rational (-2)-curves, and let
$[\eta]\in H_2(X, \Z)$ be a homology class of a section.
Denote by $[Z_1], ..., [Z_k]\in H_2(X,\Z)$ the homology classes
of the distinct exceptional (-2)-components $Z_i$ of a special fiber $F$,
and let $n_1, ..., n_k$ be a collection of integers.
Then $A:=\eta + \sum_i n_i [Z_i]$ is a homology class
whose intersection with a general fiber is 1. 
Using Thom's theorem as in Subsection \ref{_section_Nron_Subsection_},
Step 2, we can represent $A$ by a smooth 2-dimensional
submanifold of $X$. The averaging procedure,
applied to this situation, produces a section $S$
of the corresponding N\'eron model, which intersects
only one of the irreducible components of $F$,
and $A\cap [Z_i]= n_i [Z_i]^2 + \eta \cap [Z_i]$,
while $S$, being a section, might intersect only one of $Z_i$,
with multiplicity 1.

\hfill

{\bf Proof of \ref{_sections_from_CP^1_Theorem_} (ii). Step 1:}
The idea of the proof is simple. 
Let $C\subset \C P^n$ be a general line, and $M_C\subset
M$ its preimage, $M_C:= \pi^{-1}(C)$.
In \ref{_averaging_class_Theorem_}
we express the homology class of a smooth 
section $S'$ of $\pi:\; M_C \to C$, obtained
by the averaging procedure, in terms of $\alpha \in H_{2n}(M, \Z)$.
By \ref{_corollary_S'_independent_Corollary_}, the
 difference $[S]-[S']$ is expressed intrinsically
in terms of the homology, and does not depend on the choice
of the line $C\subset \C P^n$. 
We consider this class
as a section of an appropriate local system over the space
of general lines in $\C P^n$, and obtain that it is invariant
under the monodromy. Then we show that the monodromy-invariant
classes are obtained as cup products, using Deligne's theorem
on invariant cycles.

\hfill

{\bf Step 2:}
Let $B_1$ be the space of all lines in $\C P^n$ transversal
to the discriminant $D$, and ${\Bbb B}$ the
local system on $B_1$, taking value
 $H^n(M_C, \Q)$ at each line $C\in B_1$.
Using Poincar\'e duality,
we intepret the class $[S]-[S']\in H_2(M_C, \Q)$
as an element in $H^{2n}(M_C, \Q)$. This
gives a section $u$ of ${\Bbb B}$.

We use Deligne's invariant cycle theorem 
(\cite[Theorem 4.1.1]{_Deligne:Hodge_II_}), which
claims the following. Let $A_0 \to B_0$ be a 
proper holomorphic submersion of quasiprojective manifolds,
$F$ a smooth fiber, and $A$ a smooth projective compactification of $A_0$.
Then the $p$-th rational cohomology of the fibers of $\pi$
are fibers of the Gauss-Manin local system over $B_0$.
Let $\eta \in H^p(F)$ be a monodromy invariant
class in the cohomology of the fiber. Then $\eta$ belongs
to the image of the restriction map $H^p(A) \to H^p(F)$.

We apply this result to the fibration
$\Psi:\; A \to \Gr(2,n+1)$, where $\Gr(2,n+1)$ is the Grassmannian
of all lines in $\C P^n$, and 
$A=\{(m\in M, C\in \Gr(2,n+1),\ \ |\ \  \pi(m)\in C\}$,
and $A_1\subset A$ the space of all $(m, C)\in A$
such that $C\in B_1$. The local system ${\Bbb B}$
is identified with $R^{2n}\Psi_*(\Q_{A_1})$,
that is, with the Gauss-Manin local system of $H^{2n}(M_C, \Q)$
over each $C$. Since $u$ is a parallel section of ${\Bbb B}$,
its value $u\restrict C\in H^{2n}(M_C, \Q)$ at each $C$
is monodromy invariant. Therefore,
Deligne's theorem implies that $u\restrict C$ is obtained
as a restriction of a cohomology class $U \in H^{2n}(A, \Q)$.

The manifold $A$ is fibered over $M$ with the fiber
at each point $m\in M$ consisting of all lines
passing through $\pi(m)$, that is, $\Gr(2,n)$-fibered
over $M$. Therefore, 
$H^*(A)= H^*(M) \otimes H^*(\Gr(2,n))$
(\cite{_Deligne:degenerate_}).
In the next step, we are going to show that
the restriction $H^{2n}(A, \Q)\to H^{2n}(M_C, \Q)$
is factored through $H^{2n}(M, \Q)$. Therefore,
we can consider $U$ as an element in $H^{2n}(M,\Q)$.
Since the Poincar\'e dual class $([S]-[S'])^*$
is obtained as a restriction of $U$, 
this implies that the homology class $[S]-[S']$ is obtained
as a cup product of $M_C$ with the dual 
class $U^*\in H_{2n}(M,\Q)$.

\hfill

{\bf Step 3:} 
Fix $C\in B_1$.
Consider the forgetful map $\Phi:\; A \to \Gr(2, n+1)$
taking $(m, C)$ to $C$. Clearly, $\Phi^{-1}(C)=M_C$.
This identifies $M_C$ with a fiber of $\Phi$.
Using the projection $\Psi:\; A \to M$,
we consider $H^*(M)=\Psi^*(H^*(M))$ as a subspace in $H^*(A)$.
We are going to show that the image of 
the restriction map $r:\; H^*(A)\to H^*(M_C)$
is contained in the restriction of $\Psi^*H^*(M)$.
As explained in Step 2, $\im r\subset \Psi^*H^*(M)$ implies
 \ref{_sections_from_CP^1_Theorem_} (ii).

\hfill

{\bf Step 4:}
Fix a general point $x\in C$, and
let $M_x\subset A$ be the set of all pairs 
$M_x:= \{(m\in M, C'\ni \pi(m)) \ \ | \ \  C'\ni x\}$.
Clearly, $M_x$ is isomorphic to the blow-up of $M$
in $\pi^{-1}(x)$ and contains $M_C$. The restriction map $H^*(A) \to H^*(M_C)$
is factored through $M_x\supset M_C$. 

Denote by $\sigma:\; M_x \to M$ the blow-up map.
In Step 3, we deduced \ref{_sections_from_CP^1_Theorem_} (ii)
from the relation $\im r\subset \Psi^*H^*(M)$.
Let $r_1:\; H^*(M_x)\to H^*(M_C)$ be the restriction map.
Then $\im r\subset \Psi^*H^*(M)$ would follow if we show that the
following diagram is commutative:
\[
\begin{diagram}
H^*(M_x) & \rTo^{\sigma_*} & H^*(M)\\
&\rdTo^{r_1} &\dTo^{r_2}\\
& & H^*(M_C),
\end{diagram}
\]
where $r_1$ and $r_2$ are restriction maps, and 
$\sigma_*$ is pushforward in cohomology.
This diagram is clearly commutative, because
$\sigma_*$ is the fiberwise integration map,
and $M_C$ intersects each fiber of $\sigma$ 
at most once.
\endproof

\subsection{Constructing a smooth 
section over a curve (without N\'eron models)} 
\label{_curve_no_nron_Subsection_}
%%%%%%%%%%%%%%%%%%%%%%%%%%%%%%%%%%%%%%%%%%%%%%%%%%%%%%%%%%%%

We give another proof of \ref{_sections_from_CP^1_Theorem_} (i)
using the idea communicated to us (after the first draft of this
paper was published) by J\'anos Koll\'ar. This way, we 
can get rid of N\'eron models in \ref{_sections_from_CP^1_Theorem_} (i);
this also makes it unnecessary to assume that $\pi:\; \pi^{-1}(C)\to C$
has local sections.

\hfill

%%%%%%%%%%%%%%%%%%%%%%%%%%%%%%%%%%%%%%%%%%%%%%%%%%%%%%%%%%%%
\theorem\label{_sections_from_CP^1_no_Neron_Theorem_}
Let $\pi:\; M \to X$ be a hyperk\"ahler manifold equipped with a
Lagrangian fibration, with $X$ smooth (and hence,
by Hwang's theorem, isomorphic
to $\C P^n$), and $C\subset X$ a general line in $X=\C P^n$. 
Assume that the general fiber of the Lagrangian projection $\pi$ 
is primitive (that is, not divisible) in the integral homology. 
Then there exists a smooth section
$\sigma:\; C \arrow M$.

\hfill

{\bf Proof. Step 1:} Let $M_C:=\pi^{-1}(C)$. Like in the proof of
\ref{_sections_from_CP^1_Theorem_}, Step 2, we can assume that
$M_C$ contains a smooth submanifold $Z$ such that $[Z]\cap F=1$,
where $F=\pi^{-1}(z)$ for a general point $z\in C$.
Let $F_1$ be a singular fiber of $\pi$.
In the first and second step, we show that $Z$ can be chosen
in such a way that for some neighbourhood $U\supset F_1$,
$Z\cap U$ is a disconnected union of holomorphic subvarieties
(possibly with an opposite orientation).

If $F_1$ is reduced, this observation is clear: we 
deform $Z$ in such a way that it intersects $F_1$
transversally in its smooth points, and notice that
every smooth section of the projection $(z_1, z_2, ..., z_n) \mapsto z_1$
is homotopic to a holomorphic section near the point
$(0, z_2, ..., z_n)$.

{\bf Step 2:} 
If $F_1$ is not reduced, locally around a smooth point
$m$ of its reduction the projection $\pi:\; M_C \to C$
can be written as $(z_1, z_2, ..., z_n) \mapsto z_1^k$,
where $k$ is the multiplicity of this component of $F_1$.
Then, as in Step 1, we can deform $Z$ in such a way that
it is holomorphic in the coordinates  $(z_1, z_2, ..., z_n)$
near $m$. In a small neighbourhood of $\pi(m)\in C$, then,
$Z$ becomes a ramified cover, isomorphic to the graph
of the multivalued function $t\to (\sqrt[k]t, z_2, ..., z_n)$.

{\bf Step 3:} Let $Z_0$ be the intersection of $Z$
with the set of all smooth fibers of $\pi:\; M_C \to C$.
These fibers are complex tori, hence we can
apply to $Z_0$ the averaging procedure
defined in \ref{_sum_in_torsor_indep_Claim_},
obtaining a section of $\pi:\; M_C \to C$
over the set $C_0$ of all regular values of this map.
Since near the critical values $Z$ was a holomorphic
multivalued function, and the averaging procedure defines
a meromorphic map $M_C\times_C M_C\times_C ... \times_C M_C\to M_C$, we
obtain a section $S$ of $\pi:\; \pi^{-1}(C) \to C$ which is 
meromorphic near all singular fibers. Since the fibers of
$\pi$ are compact, $S$ is holomorphic at 
each critical value of $\pi:\; \pi^{-1}(C) \to C$ 
by the valuative criterion of properness. Therefore,
this section is smooth on $C$. \endproof

%%%%%%%%%%%%%%%%%%%%%%%%%%%%%%%%%%%%%%%%%%%%%%%%%%%%%%%%%%%%

\section{Holography principle for ample rational curves}

%%%%%%%%%%%%%%%%%%%%%%%%%%%%%%%%%%%%%%%%%%%%%%%%%%%%%%%%%%%%

%%%%%%%%%%%%%%%%%%%%%%%%%%%%%%%%%%%%%%%%%%%%%%%%%%%%%%%%%%%%%%%%%
\subsection{Barlet spaces: an introduction}
%%%%%%%%%%%%%%%%%%%%%%%%%%%%%%%%%%%%%%%%%%%%%%%%%%%%%%%%%%%%%%%%%

Let $Z\subset M$ be a compact subvariety in a complex
manifold, and ${\cal Z}$ the set of all complex deformations
of $Z$, that is, the set of all complex subvarieties
which occur as fibers in a flat holomorphic family of subvarieties
containing $Z$. There are two basic ways to 
equip ${\cal Z}$ with a structure of a complex variety:
one could consider the Douady space of subschemes
or the Barlet space of cycles in $M$. It won't make
much difference for the present purposes, but we will
focus on the Barlet spaces.
For an introduction, reference and more details
on Barlet spaces, see 
\cite{_Magnusson:cycle_,_Barlet_Magnusson:book_,KL}. 

\hfill

Equip $M$ with a Riemannian metric; for the sequel, 
it does not matter which one we chose. Let ${\goth S}$ be the set
of all closed subsets of $M$. The  {\bf Hausdorff
topology} on ${\goth S}$ is defined as follows.
Fix $\epsilon >0$ and a closed subset $Z\subset M$. The base of 
the Hausdorff topology is the set $U_\epsilon (Z)$
of all $Z'\in {\goth S}$ such that the
$\epsilon$-neighbourhood of $Z$ contains $Z'$, and the $\epsilon$-neighbourhood of $Z'$
contains $Z$. It is not hard to see that Hausdorff topology
does not depend on the choice of a Riemannian metric.

\hfill

A complex variety is called {\bf equidimensional}
if all its irreducible components have the same dimension.
A {\bf subvariety with multiplicities}, or {\bf a complex
analytic cycle} is an equidimensional subvariety $Z\subset M$
with ``mutiplicities'' (positive integer numbers)
attached to each irreducible component.
We represent complex analytic cycles by a sum $\sum_i n_i Z_i$,
where $Z_i$ are its irreducible components, and $n_i$ the multiplicities.
For every complex analytic cycle, its {\bf homology class}
is $\sum_i n_i [Z_i]$, where $Z_i\subset M$ is a complex
subvariety.

\hfill

Fix an integer homology class $\eta \in H_{2i}(M,\Z)$.
The Barlet space parametrizes complex analytic cycles
in $M$ with a fixed homology class. The Barlet space
is equipped with the topology which is induced
by the Hausdorff topology on the set of subsets.
As shown in e.g. \cite{_Barlet_Magnusson:book_}, the Barlet space is equipped
with a complex structure, which is compatible
with the Hausdorff topology.

\hfill

\definition
The {\bf Douady space}
is the moduli of equidimensional dimensional complex subspaces of $M$ (``subspaces''
here means that we allow nilpotents in the structure sheaf).
The {\bf Hilbert-Chow morphism} is the forgetful map
$F:\; {\cal D} \arrow {\cal B}$ from the Douady space to the Barlet space, with
multiplicity of each component of $F(\tilde Z)$ given by the
rank of the structure sheaf of $Z$ in the general point of this component.

\hfill

\example
For $\dim Z=0$, the Douady space is the Hilbert space
$M^{[n]}$ of 0-dimensional subschemes, and the Hilbert-Chow morphism 
is the natural projection $M^{[n]}\arrow \Sym^n M$.

\hfill

For the present purposes, the following theorem
is most important. Fix a Hermitian form $\omega$ on $M$.
Recall that the volume $\Vol_\omega(Z):= \int_Z \omega^{\dim_\C Z}$
is well defined for any complex subvariety, even singular,
finite when $Z$ is compact,
and equal to the cohomology pairing $\langle [Z], \omega^{\dim_\C Z}\rangle$
when $\omega$ is closed (\cite{_Stolzenberg_}). 
For each complex analytic cycle $Z=\sum_i n_i Z_i$
of dimension $k$, we define its {\bf Hermitian volume}
as $\Vol_\omega(Z)=\sum_i n_i \int_{Z_i} \omega^k$.
It is easy to see that the function $\Vol_\omega$
is continuous on the Barlet space.

\hfill

\theorem\label{_Bishop_compactness_Theorem_}
 (E. Bishop's compactness theorem)\\
Let $(M,\omega)$ be a compact Hermitian manifold,
and $B_\eta$ the Barlet space of all cycles homologous
to $\eta\in H_{2k}(M)$. Denote by $B_\eta(C)$
the subset of $B_\eta$ consisting of all
cycles $Z$ with $\Vol_\omega(Z)\leq C$.
Then $B_\eta(C)$ is compact.

\proof \cite{_Magnusson:cycle_,_Barlet_Magnusson:book_}. 
\endproof

\hfill

\remark
This theorem also holds for pseudoholomorphic
curves on compact Hermitian almost complex manifold;
in this context it is known as {\em Gromov's compactness
theorem} (\cite{_Gromov:curves_}).

%%%%%%%%%%%%%%%%%%%%%%%%%%%%%%%%%%%%%%%%%%%%%%%%%%%%%%%%%%%%%%%%%
\subsection{Rationally connected varieties}
%%%%%%%%%%%%%%%%%%%%%%%%%%%%%%%%%%%%%%%%%%%%%%%%%%%%%%%%%%%%%%%%%

We give a small introduction to rational connectedness, 
used in the next subsection. For a proper introduction
to this subject, see \cite{_Campana:RC_,_KMM:RC_}.

\hfill

\definition
A compact complex manifold $M$ is called {\bf rationally connected}
if any general $x,y\in M$ can be connected by a sequence
of intersecting rational curves, and all irreducible
components of Barlet space of curves on $M$ are compact.

\hfill

\remark
The condition ``all irreducible
components of Barlet space of curves on $M$ are compact''
would follow, for instance, if $M$ is Fujiki class C,
that is, bimeromorphic to K\"ahler (\cite{_Fujiki:Douady_82_}).
However, this condition is more general (\cite{_V:curves_twistor_}).

\hfill

\definition
Recall that a smooth rational curve in a complex manifold
is called {\bf ample} if its normal bundle is $\bigoplus_i \calo(k_i)$,
with all $k_i$ strictly positive.

\hfill

%%%%%%%%%%%%%%%%%%%%%%%%%%%%%%%%%%%%%%%%%%%%%%%%%%%%%%%%%%%%
\claim\label{_rationally_conne_Claim_}
Let $M$ be a compact K\"ahler manifold containing 
an ample rational curve. Then any $N$ points $z_1, ..., z_N$ can
be connected by an ample rational curve.

{\bf Proof:} This fact is well known in algebraic geometry
(see \cite{_Kollar:curves_}). However, its proof is valid
for all K\"ahler manifolds. \endproof

\hfill

%%%%%%%%%%%%%%%%%%%%%%%%%%%%%%%%%%%%%%%%%%%%%
\claim\label{_ratl_conne_then_ample_Claim_}
A compact K\"ahler 
manifold which contains an ample rational curve
is rationally connected.

\proof See e.g. \cite[Lemma 1.2]{_Harris:RC_}
\endproof

\hfill

%\definition
%Let $M$ be a compact K\"ahler manifold, or a variety which
%is bimeromorphic to a compact K\"ahler manifold.
%{\bf The MRC quotient} of $M$ is a dominant rational
%map $M \dashrightarrow Z$ to a compact variety such that its general
%fibers are bimeromorphic to rationally connected
%manifolds, and general points of $Z$ are not contained
%in any rational curve.
%
%\hfill
%
%As shown in \cite{_Campana:RC_,_KMM:RC_},
%the MRC quotient exists, and is unique up to 
%bimeromorphic equivalence.

Clearly, 
rational connectedness is invariant under bimeromorphic
maps (\cite[Chapter IV, Proposition 3.3.3]{_Kollar:curves_}). 
However, it is not that easy to define rational
connectedness for singular varieties. For example,
consider a rank 2 holomorphic vector bundle $B$ over $Y$, 
where $Y$ is a smooth compact complex curve of genus $\geq 1$.
Let $X:= {\Bbb P}B$ be its projectivization. It is easy to check
that any rational curve in $X$ is projected to a point in $Y$,
hence $X$ is not rationally connected. Now, assume that
$B=L \oplus \calo_Y$, where $L^*$ is ample. Then
the zero section of the projection $X \to Y$
can be blown down, resulting in a variety $X_1$
which can be also obtained as a compactification 
of the affine cone $C(Y)$. The rational curves
in $X$ now all intersect in the singular point
of the origin of the cone, but the resolution of $X_1$ is not
rationally connected.

In the singular setting, definition of rational connectedness
can be given as follows. A variety is {\bf rationally connected}
if its resolution of singularities is rationally connected.
By \cite[Chapter IV, Proposition 3.3.3]{_Kollar:curves_},
this is equivalent to the following definition.

\hfill

%%%%%%%%%%%%%%%%%%%%%%%%%%%%%%%%%%%%%%%%%%%%%%%%%%%%%%%%%%
\definition\label{_singular_RC_Definition_}
Let $M$ be an irreducible compact complex variety $M$ 
bimeromorphic to a compact K\"ahler manifold.
Then $M$ is called {\bf rationally (chain) connected}
if for any proper subvariety $E\subset M$,  
any two points in $M$ can be connected by a sequence
of intersecting rational curves which intersect in points of 
$M \backslash E$.

\hfill

However, for our purposes, an easier criterion is sufficient.

\hfill

%%%%%%%%%%%%%%%%%%%%%%%%%%%%%%%%%%%%%%%%%%
\claim\label{_RC_one_curve_Claim_}
Let $M$ be a compact complex variety, bimeromorphic to a 
compact K\"ahler manifold, such that two generic 
points $x, y\in M$ can be connected by 
a smooth rational curve $S_{x,y}$. Then 
$M$ is rationally connected.

\hfill

\proof Let $M_1\stackrel \pi\arrow M$ be the K\"ahler resolution of singularities.
For $x, y\in M_1$ general, the curve $S_{\pi(x),\pi(y)}$ 
does not belong to the exceptional locus of $\pi$,
hence it can be lifted to $M_1$. Therefore,
$M_1$ is rationally connected.
\endproof

%%%%%%%%%%%%%%%%%%%%%%%%%%%%%%%%%%%%%%%%%%%%%%%%%%%%%%%%%%%%%%%%%
\subsection{Holography principle with K\"ahler target space}
%%%%%%%%%%%%%%%%%%%%%%%%%%%%%%%%%%%%%%%%%%%%%%%%%%%%%%%%%%%%%%%%%

To proceed further, we need to generalize the result of
\cite{_V:holography_}, extending the ``holography principle''
to a K\"ahler target space. 
This is the main result of the present section.

\hfill

%%%%%%%%%%%%%%%%%%%%%%%%%%%%%%%%%%%%%%%%%%%%%%%
\theorem \label{_holo_for_Kahler_Theorem_}
Let $S\subset M$ be an ample rational curve in a 
$M=\C P^n$, and $U\supset S$ a connected
neighbourhood of $S$. 
Consider an open, connected neighbourhood $U_S\subset M$ of $S$.
Then any holomorphic map $\phi:\; U_S \arrow X$ from $U_S$ to a K\"ahler manifold
$X$ can be extended to a meromorphic map $M \dashrightarrow X$.

\hfill

\remark
For \ref{_holo_for_Kahler_Theorem_},
K\"ahlerness of $X$ is sufficient, but not necessary.
As seen from the proof, for \ref{_holo_for_Kahler_Theorem_}
to hold it suffices to establish that
all components of the Barlet spaces of curves on $X$ and $M$
are compact. This holds true if $X$ admits a Hermitian
form $\omega$ which is {\em pluriclosed} (i.e., it satisfies $dd^c\omega=0$),
or such that $-dd^c \omega$ is a positive form 
(\cite{_V:curves_twistor_}). 

\hfill

\remark 
It is not difficult to find a counterexample
to \ref{_holo_for_Kahler_Theorem_} when $X$
is not simply connected and non-K\"ahler. 
For instance, let $X=\Tw(Y)$ be the twistor space
for a 2-dimensional complex torus $Y$ (see e.g. 
\cite{_V:curves_twistor_}). This is an example
of a complex manifold with Barlet components
of rational curves non-compact (indeed, in
\cite{_V:curves_twistor_} it is shown that
all of these components are Stein).
The universal cover
of $\Tw(Y)$ is the twistor space for $\C^2={\Bbb H}$,
which can be identified with the total space
of $\calo(1)^2$ (\cite{_JV:Instantons_}).
Since $\Tw(S^4)$ is $\C P^3$ and $\C^2$
is conformally equivalent to $S^4$ without a point,
the space $\Tw(\C^2)$ is biholomorphic
to $\C P^3\backslash \C P^1$ (see 
\cite{_JV:Instantons_} for another argument 
leading to this conclusion). Then 
an open neighbourhood of $\C P^1\subset \C P^3$
can be embedded into $\Tw(Y)$. However, if this
embedding can be extended to the whole $\C P^3$,
the space $\Tw(Y)$ would be Moishezon, being
an image of a proper holomorphic 
map from a Moishezon space (see e.g. 
\cite{_Fujiki:Douady_82_}). This is
impossible, because the algebraic
dimension of any twistor space of
a compact hyperk\"ahler manifold is one 
(\cite{_V:holography_}).
 
\hfill

We expect that a similar
counterexample would exist when $X$ is non-K\"ahler but
simply connected. 

\hfill

When $X$ is Moishezon, \ref{_holo_for_Kahler_Theorem_}
is implied by the following theorem (see also 
\cite[Theorem 3.4]{_V:holography_}; if we deal with formal
schemes instead of complex manifolds, this result was 
established in \cite{_Hartshorne:cohomo-1968_} and in 
\cite{_Hironaka_Matsumura:formal_,Hironaka}). The following
result requires the manifold $M$ to be homogeneous and rational; for
non-homogeneous case, see \cite[Theorem 3.4]{_V:holography_}.

\hfill

Let $Z\subset M=G/H$ be a subvariety in a complex homogeneous manifold, where $G$
is an algebraic Lie group, 
$p \in Z$ a point, and $G_{p,Z}$ the set of all elements in $G$ which carry
$p$ to a point in $Z$. Denote by $G_Z$ the group generated by the set $G_{p,Z}$;
clearly, $G_Z$ is an algebraic subgroup of $G$, independent from the choice
of $p\in Z$. We say that $Z$ {\bf generates M} if $G_Z=G$.

\hfill
 
%%%%%%%%%%%%%%%%%%%%%%%%%%%%%%%%%%%%%%%%%%%%%%%%
%\theorem \label{_holo_restrictions_Theorem_}
%Let $S\subset M$ be an ample
%rational curve in a simply connected complex manifold,
%which is covered by smooth, ample deformations of $S$.
%Consider a neighbourhood $U\supset S$ which is contained
%in a union ${\goth S}$ of all deformations of $S$
%intersecting $S$.\footnote{This assumption holds for all $U$ 
%if $M$ is compact and K\"ahler, 
%or if all components of the Barlet spaces of curves for $M$ are compact
%\cite[Remark 3.2]{_V:holography_}.}
%Then for any smaller 
%open neighbourhood $V\subset U$ of $S$,
%the restriction map $\Mer(U) \arrow \Mer(V)$ 
%is an isomorphism.

%%%%%%%%%%%%%%%%%%%%%%%%%%%%%%%%%%%%%%%%%%%%%%%%%%%%%%%%%%%%
\theorem\label{_Chow_rational_extension_Theorem_}
Let $M$ be a rational homogeneous algebraic variety, $Z$ a
complete algebraic subvariety of positive dimension in $M$
which generates $M$, and let $U$ be a connected neighborhood
of $Z\subset M$. Then any meromorphic map from $U$ to 
an algebraic variety $X$ is a restriction of a rational
map from $M$ to $X$. 

\proof This is \cite[Theorem 2, Theorem 3]{_Chow:meromorphic_}.
Note that ``algebraic map'' in \cite[Theorem 2]{_Chow:meromorphic_}
refers to multi-valued maps given by an algebraic correspondence,
which is shown to be rational in \cite[Theorem 3]{_Chow:meromorphic_}.
\endproof

\hfill

{\bf Proof of \ref{_holo_for_Kahler_Theorem_}:}
The main idea of the proof of \ref{_holo_for_Kahler_Theorem_}
is the following observation. Consider an irreducible
complex subvariety $B_0$ in the Barlet space of curves
on $X$ containing $n\phi(S')$, where $n$
is a positive integer (representing multiplicity)
and  $S'\subset U_S$ a general deformation of $S$ in $U_S$.
The corresponding {\bf Barlet locus} $L_0$ is defined as the union of all
curves $C\in B_0$. Assume that $L_0$ is rationally connected.
Then $L_0$ is Moishezon (\cite{_Campana:rational_lines_}) 
and its normalization is simply connected (\cite[Corollary 2.32]{_Debarre:RCV_};
see also \cite{_Campana:RC_,_KMM:RC_}),
hence \ref{_Chow_rational_extension_Theorem_} can be applied,
showing that the map $\phi:\; U_S \arrow L_0$
can be extended to a meromorphic map $M \dashrightarrow L_0\subset X$.

\hfill

{\bf Step 1:}
We take $n=3$. Let $B_0$ be the set of all deformations
of $3\phi(S)$ (that is, $\phi(S)$ with multiplicity 3)
which intersect $\phi(S)$ and intersect $\phi(U_S)$
in an open set. The first condition is closed.
To get a hold on the second condition, 
denote by $s$ the intersection of a curve
$S_1\in B_0$ and $\phi(S)$.
Clearly, $S_1\cap \phi(U_S)$ is open 
if and only if $S_1$ is tangent 
of infinite order\footnote{That is, its $k$-jet in $s$ 
belongs to the set of $k$-jets of curves in $U_S$, for all $k$.}
to $U_S$ in $s$. This condition
is closed, because it is obtained by
taking an intersection of closed subvarieties enumerated by 
the order of tangency $k=1, 2, 3, \dots$. Therefore, $B_0$
is closed in the Barlet component $B$ of all
deformations of $3\phi(S)$. 
From Bishop's compactness theorem
it follows that $B_0$ is compact. 
From Remmert's proper
map theorem, we obtain that 
the corresponding Barlet locus
$L_0\subset X$ is complex analytic. 
To finish the proof, it remains
to show that $L_0$ is rationally connected.
We will prove that two generic points of $L_0$
are connected by a smooth curve which belongs
to $B_0$; by \ref{_RC_one_curve_Claim_}, this implies
the rational connectedness of $L_0$.

\hfill

{\bf Step 2:}
Since $S$ is ample, its normal bundle $NS$
is a direct sum of $\calo(i)$ with $i\geq 1$.
Therefore, for any two distinct points $x, y\in S$,
there exists neighbourhoods $U_x\ni x, U_y\ni y$ in $U_S$
such that any two points $a\in U_x$, $b\in U_y$
can be connected by a deformation of $S$
which belongs to $U_S$. When $NS=\calo(1)^n$,
this curve is unique, hence it cannot intersect $S$
when $a, b$ are generic.

A union of two deformations $S_1\ni a, S_2\ni b$ of $S$ which 
intersect can be deformed to a smooth deformation $S_3$
of $2S$ by the standard technique (\cite{_Kollar:curves_}).
Moreover, we can perform this deformational smoothing
in such a way that $S_3$ contains $a$ and $b$.
This is the idea used to prove \ref{_rationally_conne_Claim_}.

Start from $a, b\in U_S$ generic and sufficiently close to $S$, 
and let $S_1$ connect $a$ to $x$,
$S_2$ connect $b$ to $x$, and $S_3$ connect $z\in S$ to $t\in S_1$.
Smoothing the union $S_1\cap S_2\cap S_3$ as above, we obtain
a smooth curve in $U_S$ passing through $a, b$ and 
$z\in S$. By construction, this curve belongs to the Barlet space
of deformations of $3 S$.

\hfill

{\bf Step 3:}
Let $L_0\subset X$ be the subvariety constructed in Step 1.
By construction, $\phi(U_S)$ is an open neighbourhood of 
$S\subset L_0$. Therefore, to prove that $L_0$ is rationally
connected, it would suffice to find a smooth rational curve
in $L_0$ connecting two general points in $\phi(U_S)$.
This was done in Step 2. As shown in Step 1,
this implies that $L_0$ is rationally connected,
and finishes the proof of
\ref{_holo_for_Kahler_Theorem_}.
\endproof

%%%%%%%%%%%%%%%%%%%%%%%%%%%%%%%%%%%%%%%%%%%%%%%%%%%%%%%%%%%%

\section{Meromorphic sections via Dolbeault currents} 
\label{_Meromorphic_Section_}

%%%%%%%%%%%%%%%%%%%%%%%%%%%%%%%%%%%%%%%%%%%%%%%%%%%%%%%%%%%%

%%%%%%%%%%%%%%%%%%%%%%%%%%%%%%%%%%%%%%%%%%%%%%%%%%%%%%%%%%%%
\subsection{The Dolbeault class of a section} 
\label{_Dolbeault_classes_Subsection_}
%%%%%%%%%%%%%%%%%%%%%%%%%%%%%%%%%%%%%%%%%%%%%%%%%%%%%%%%%%%%

\definition
Let $\pi:\; M \arrow X$ be a Lagrangian fibration,
and $\sigma:\; X \arrow M$ a section.
By \cite[Proposition 2.1]{_BDV:sections_}
(see also \ref{_restriction_2,0+1,1_Proposition_} below), 
 $\sigma^*\Omega$ is a $(2,0)+(1,1)$ closed
form on $X$. Assume that $H^2(X)$ is generated by
a closed $(1,1)$-form $\omega$.
Then $\sigma^*\Omega-\lambda \omega$
is exact for some $\lambda\in \C$.
Let $\alpha$ be a 1-form which satisfies 
$d\alpha=\sigma^*\Omega-\lambda \omega$.
Then $\bar\6\alpha^{0,1}=0$.
The {\bf Dolbeaut class} of $\sigma$ 
is the class of $\alpha^{0,1}$ in
$H^1(\Lambda^{0,*}X,\bar\6)=H^1(X, \calo_X)$.

\hfill

\proposition \label{_Dolbealt_vanishes_then_extend_Proposition_}
Let $\pi:\; M \arrow \C P^n$
be a Lagrangian fibration, 
$C\subset \C P^n$ a smooth curve,
$U_C$ its neighbourhood and 
$\sigma_U:\; U_C\arrow M$
a section. Assume that the Dolbeault class 
$[\alpha^{0,1}]\in H^1(U_C, \calo_{U_C})$ of $\sigma_U$ vanishes.
Then $\sigma^*_U \Omega$ restricted to a smaller
  neighbourhood of $C$ can be extended to a closed 
$(2,0)+(1,1)$-form on $\C P^n$.

\hfill

\proof
Let $\alpha$ be a 1-form which satisfies 
$d\alpha=\sigma^*\Omega-\lambda \omega$, where $\omega$ is
the Fubini-Study form.
Since the $\bar\6$-class of $\alpha^{0,1}$ vanishes,
there exists $f\in C^\infty M$ such that 
$\alpha':=\alpha-\bar\6 f$ is of type $(1,0)$.
Shrinking $U_C$ if necessary,
we can extend $f$ to a function $\tilde f$ on $\C P^n$
and $\alpha'$ to a $(1,0)$-form $\tilde \alpha'$ on $\C P^n$.
Set  $\tilde \alpha:=\tilde \alpha'+\bar\6 \tilde f$.
Then $d\tilde \alpha+\lambda \omega$ is a 
closed $(2,0)+(1,1)$-form on $\C P^n$ extending $\sigma^* \Omega$.
\endproof

\hfill

\definition
Let $U\subset \C P^n$ be an open subset such that 
$H^2(U)$ is generated by the Fubini-Study form.
Let $\rho$ be a closed $(2,0)+(1,1)$-form 
on $U$, and $\rho= d\alpha + \lambda \omega$.
Then the {\bf $(0,1)$-Dolbeault class} of $\rho$ is
$[\alpha^{0,1}]\in H^{0,1}(M)$. 

\hfill

\remark
By definition, the Dolbeault class of $\sigma_U$ 
is equal to the (0,1)-Dolbeault class of $\sigma_U^*\Omega$. 

\hfill

In the conclusion of this section
(Subsection \ref{_vanishing_Subsection_}), we 
shall prove the following theorem:

\hfill

%%%%%%%%%%%%%%%%%%%%%%%%%%%%%%%%%%%%%%%%%%%%%%%%%
\theorem
Let $\pi:\; M \arrow \C P^n$
be a Lagrangian fibration, 
$C\subset \C P^n$ a general line, 
let $U_C$ be a neighbourhood of $C$, and 
$\sigma_U:\; U_C\arrow M$
 a smooth section which satisfies 
\ref{_sections_from_CP^1_Theorem_} (ii). Then the
Dolbeault class of $\sigma_U$ vanishes.

\proof \ref{_Dolbeault_vanishes_Theorem_}. \endproof

\hfill

This theorem completes the last step of 
the proof of \ref{_Sections_main_Theorem_},
the main result of this paper.
Before we prove this theorem, we need
to make a digression to the land of currents.

%%%%%%%%%%%%%%%%%%%%%%%%%%%%%%%%%%%%%%%%%%%%%%%%%%%%%%
\subsection{Currents on complex manifolds}
%%%%%%%%%%%%%%%%%%%%%%%%%%%%%%%%%%%%%%%%%%%%%%%%%%%%%%

For a rich and detailed introduction to currents and their use in 
algebraic geometry, see \cite{_Demailly:analytic_}.
We give a few preliminary definitions and cite basic results to facilitate
the proof of vanishing of the Dolbeault class
of a section, introduced earlier in Subsection \ref{_Dolbeault_classes_Subsection_}.
For all missing proofs, please see \cite{_Demailly:analytic_}.

\hfill

\definition 
Let $F$ be a Hermitian bundle with connection $\nabla$, 
on a Riemannian manifold $M$ with Levi-Civita connection, and 
\[
\|f\|_{C^k}:=\sup_{x\in M}\left(|f| + |\nabla f| + ... + |\nabla^kf|\right)
\]
the corresponding $C^k$-norm 
defined on smooth sections with compact support.

\hfill

\remark
The $C^k$-topology is independent from the choice of a
connection and a metric (\cite{_Hamilton:Moser_}).

\hfill

\definition
A {\bf generalized function} is a functional on the top forms
with compact support, which is continuous in one of the $C^i$-topologies.

\hfill

\definition
A {\bf $k$-current} is a functional on the $(\dim_\R M-k)$-forms
with compact support, which is continuous in one of the $C^i$-topologies.

\hfill

\remark Currents are forms with coefficients in generalized functions.

\hfill

\definition
The space of currents is equipped with the {\bf weak topology} 
(a sequence of currents converges if it converges on all
forms with compact support). The space of currents with 
this topology is a {\bf Montel space} (barrelled,
locally convex, all bounded subsets are precompact). 
Montel spaces are {\bf reflexive} (the map
to its double dual with the strong topology is an isomorphism).

\hfill

\definition
On a complex manifold, {\bf $(p,q)$-currents} are $(p,q)$-forms with 
coefficients in generalized functions.

\hfill

\remark In the literature, this is sometimes called
``$(n-p,n-q)$-currents''; we follow the conventions from
\cite{_Demailly:analytic_}.

\hfill

\claim De Rham differential is continuous on currents, and 
the Poincar\'e lemma holds for currents, as well as the
Poincar\'e-Dolbeault-Grothendieck lemma. Hence, 
the cohomology of currents are the same
as cohomology of smooth forms,
and the $\bar\6$-cohomology groups are the same as for forms. 

\hfill

\definition
Let $f:\; X \arrow Y$ be a proper holomorphic map of complex manifolds,
$\dim_\C X = \dim_\C Y +k$, and $\alpha$ a $(p,q)$-current on $X$.
Define the {\bf pushforward} $f_* \alpha$ using
$\langle f_*\alpha, \tau\rangle := \langle \alpha, f^*\tau\rangle$,
where $\tau$ is any test-form (a smooth differential form with compact 
support). Then $f_*\alpha$ is a
$(p-k,q-k)$-current. One should think of $f_*$ as 
fiberwise integration. 

\hfill

\remark Clearly, $d f_* \alpha= f_* d\alpha$ for any smooth map,
$\6 f_* \alpha= f_* \6\alpha$ for a holomorphic map, and so on. 

\hfill

\remark Pullbacks of currents are (generally speaking) not well-defined.

\hfill 

Using currents, one can give a new proof of 
the following result of \cite{_BDV:sections_}.
 
\hfill 

\proposition \label{_restriction_2,0+1,1_Proposition_}
Let  $(M, I, \Omega)$ be a holomorphically symplectic manifold
equipped with a holomorphic Lagrangian fibration
$\pi:\; M \arrow X$, and let $\sigma:\; X \arrow M$ be a smooth 
section. Then the pullback $\sigma^*(\Omega)$ is a form
of type $(2,0)+(1,1)$ on $X$. 

\hfill 

\proof
Let $Z:=\im \sigma$, and let $[Z]$ be its integration current,
defined by $\langle [Z], \tau\rangle= \int_Z \tau$.
Clearly, $\sigma^*(\Omega)= \pi_*(\Omega\wedge [Z])$.
Let $\dim_\C M=2n$. To prove that $\pi_*(\Omega\wedge [Z])$
is of type $(2,0)+(1,1)$ it suffices to show that 
$\pi_*(\Omega\wedge [Z])\wedge \tau=0$ for any $(n, n-2)$-form $\tau$.
This would follow if we prove that $\pi^* \tau\wedge \Omega=0$.
This assertion is clear, because in coordinates 
$\Omega= \sum dp_i \wedge dq_i$, where $dp_i$ are 
pullbacks of the coordinates on $X$, and 
$\tau$ is proportional to $dp_1\wedge dp_2 \wedge ... \wedge dp_n$.
\endproof

%%%%%%%%%%%%%%%%%%%%%%%%%%%%%%%%%%%%%%%%%%%%%%%%%%%%%%%%%%%%

\subsection{Vanishing of the cohomology class of the Dolbeault current}
\label{_vanishing_Subsection_}

%%%%%%%%%%%%%%%%%%%%%%%%%%%%%%%%%%%%%%%%%%%%%%%%%%%%%%%%%%%%

\theorem \label{_Dolbeault_vanishes_Theorem_}
Let $\pi:\; M \arrow \C P^n$
be a Lagrangian fibration, 
$C\subset \C P^n$ a general line, 
$U_C$ a neighbourhood of $C$, and 
$\sigma_U:\; U_C\arrow M$ 
a smooth section which satisfies 
\ref{_sections_from_CP^1_Theorem_} (ii).  Then the
Dolbeault class of $\sigma_U$ vanishes.

\hfill

\pstep
Let $\pi:\; M \arrow X$ be a 
Lagrangian fibration, with  $\Omega-\pi^* \omega=d\alpha$ 
exact and $\omega$ a closed $(1,1)$-form on $X$. 
We generalize the Dolbeault class 
to any piecewise smooth multiple section 
$Z\subset M$, which is proper and 
$k$-sheeted (if we count sheets with 
orientations) over $X$. This is done as 
follows: we split $Z$ onto smooth pieces, 
restrict $\alpha$ to each, transport to 
$X$ with the appropriate sign and sum together, 
and then take the $\bar\6$-class of the (0,1)-part of resulting current.
This is equivalent to taking the current 
$\pi_*(\alpha \wedge [Z])$, using \ref{_restriction_2,0+1,1_Proposition_}
to show that its (0,1)-part is $\bar\6$-closed, 
and taking its $\bar\6$-class. 

\hfill

{\bf Step 2:} With this definition, 
we can extend the notion of a Dolbeault class 
to any multisection of a Lagrangian fibration 
$\pi:\; M \arrow \C P^n$, and, indeed, to any 
piecewise smooth homology chain $u\in H^{2n}(M, \Z)$, 
where  $\dim_\C M=2n$.  
Defined this way, the Dolbeault class
is a homology invariant: for an exact $[Z]= ds$, 
the pushforward $\pi_*(\Omega \wedge [Z])$
satisfies 
\[\pi_*(\Omega \wedge ds)= d(\pi_*(\Omega \wedge s))-
\pi_*(d\Omega \wedge s)= d(\pi_*(\Omega \wedge s)).
\]
Then the Dolbeault class of $\sigma$ is 
equal to the $(0,1)$-Dolbeault class of the current 
$d(\pi_*(\Omega \wedge s))$, but 
the (0,1)-Dolbeault class of an exact current vanishes. 

\hfill

{\bf Step 3:} Since the section $\sigma_C$ satisfies 
\ref{_sections_from_CP^1_Theorem_} (ii), its 
homology class is obtained as a cap product of $\pi^{-1}(C)$
with an integer homology class $\alpha \in H_{2n}(M,\Z)$. This implies 
that the Dolbeault class of $\sigma_U$ vanishes; indeed, a restriction 
of any class in $H^{0,1}(\C P^n)$ to $U_C$ vanishes, because $H^{0,1}(\C P^n)=0$. 
\endproof

\hfill

Combining this theorem with \ref{_holo_for_Kahler_Theorem_},
we obtain the main result of this paper. 

\hfill

%%%%%%%%%%%%%%%%%%%%%%%%%%%%%%%%%%%%%%%%%%%%%%%%%%%%%%%%%%%%
\theorem\label{_main_res_last_section_Theorem_}
Let $\pi:\; M \arrow \C P^n$ be a Lagrangian fibration on
a hyperk\"ahler manifold. Assume 
 the general fiber of the Lagrangian projection $\pi$ 
is primitive (that is, not divisible) in integer homology.
Then there exists a
genenerate twistor deformation $I_\eta$ and
a complex subvariety $S \subset (M, I_\eta)$
such that the map $\pi:\; S \arrow \C P^n$
is birational.

\hfill

\pstep
In \cite{_Soldatenkov_Verbitsky:degenerate_Kahler_}, 
it was shown that any degenerate
twistor deformation of a hyperk\"ahler manifold is K\"ahler.
Let $C\subset \C P^n$ be a general line in the base, 
and $\sigma:\; C \arrow M$ a smooth section 
constructed in \ref{_sections_from_CP^1_Theorem_}. 
We extend it to an open neighbourhood $U_C$ of $C$, 
obtaining a map $\sigma_U:\; U_C \arrow M_U$, 
where $M_U:= \pi^{-1}(U_C)$. 
We apply \ref{_dege_twistor_defo_Theorem_} to construct a degenerate 
twistor deformation $M'_U$ of $M_U$ such that the section 
$\sigma_U:\; U_C \arrow M_U'$ is holomorphic. If this degenerate twistor
deformation can be extended from $M_U'\arrow U_C$
to $M' \arrow \C P^n$, then
\ref{_holo_for_Kahler_Theorem_}, applied to the 
holomorphic map $U_C \arrow M'$, gives a rational section of $\pi$. 

\hfill

{\bf Step 2:} 
Let $\eta_U:= \sigma_U^*(\Omega)$.
By \ref{_restriction_2,0+1,1_Proposition_}, $\eta_U$ is of type $(2,0)+(1,1)$. 
To obtain $M'$, we have to extend $\eta_U$ 
to a closed form of type $(2,0)+(1,1)$ on $\C P^n$. 
As seen in Subsection \ref{_Dolbeault_classes_Subsection_}, there is a cohomological obstruction
to such an extension, which belongs to the group 
$H^1(U_C, \calo)$.\footnote{It is possible to see that
this group is infinite-dimensional.}

\hfill

{\bf Step 3:} In \ref{_Dolbeault_vanishes_Theorem_} 
we proved that 
the obstruction class in $H^1(U_C, \calo)$ (i.e., ``the 
Dolbeault class of the section $\sigma_U$'') vanishes.
Now, \ref{_Dolbealt_vanishes_then_extend_Proposition_} implies
that the degenerate twistor deformation $M_U'$ can be extended
to $M'\arrow \C P^n$.
\endproof

%%%%%%%%%%%%%%%%%%%%%%%%%%%%%%%%%%%%%%%%%%%%%%%%%%%%%%%%%%%%

\section{Appendix. Topology
of Lagrangian fibrations}
\label{_HO_Appendix_Section_}

%%%%%%%%%%%%%%%%%%%%%%%%%%%%%%%%%%%%%%%%%%%%%%%%%%%%%%%%%%%%

Throughout most of this appendix, we 
recall the results (by now, classical)
of J.-M. Hwang and K. Oguiso about geometry
of characteristic lines on the general singular
fibers of the Lagrangian fibration. In conclusion,
we use these results to compute the rational
cohomology of those fibers.

%%%%%%%%%%%%%%%%%%%%%%%%%%%%%%%%%%%%%%%%%%%%%%%%%%%%%%%%%%%%
\subsection{Classification of singular fibers 
in the work of Hwang and Oguiso} 
\label{_Hwang_Oguiso_Subsection_}
%%%%%%%%%%%%%%%%%%%%%%%%%%%%%%%%%%%%%%%%%%%%%%%%%%%%%%%%%%%%

In \cite[Theorem 1.1]{_Hwang_Oguiso:multiple_}, 
J.-M. Hwang and K. Oguiso
classified the general singular fibers of 
proper holomorphically Lagrangian fibrations.
We use their result to describe the geometry
of irreducible components of fibers in codimension 1.
To state \cite[Theorem 1.1]{_Hwang_Oguiso:multiple_}, 
we need to give some preliminary definitions,
also taken from \cite{_Hwang_Oguiso:multiple_}.

\hfill

%%%%%%%%%%%%%%%%%%%%%%%%%%%%%%%%%%%%%%%%%%%%%%%%%%%%%%%%%%%%
\definition\label{_char_leaf_Definition_}
Let $V$ be a reduction of a fiber $\pi^{-1}(z)$ 
of a Lagrangian fibration, where $z$ is a general point
in the discriminant, $\hat V$ is its normalization,
and $\tau:\; \hat V \to \Alb(V)$ the Albanese map.
As shown \cite[Theorem 1.3]{_Hwang_Oguiso:characteristic_}, 
for each connected component $\hat V_i$ of $\hat V$, the
map $\tau:\; \hat V_i\to \Alb(\hat V_i)$ is a smooth fibration,
with all fibers isomorphic to either $\C P^1$ or an
elliptic curve. The image of each of these curves in 
$V$ is called {\bf a characteristic leaf}. 
A {\bf characteristic curve} is a maximal connected
union of characteristic leaves. Given a characteristic
curve $C=\bigcup C_i$, where $C_i$ are its irreducible
components, we define the corresponding {\bf charateristic
cycle} ${\goth C}:=\Sigma_i n_i [C_i]$, where 
$n_i$ denotes the multiplicity of the irreducible
component of $V$ which contains $C_i$.

\hfill

In \ref{_HO_classification_Theorem_} below, we repeat
\cite[Theorem 1.1]{_Hwang_Oguiso:multiple_}
adding the description of singular fibers
of complex surfaces taken from \cite[\S V.7]{_BHPV_}.

In his classification of special fibers of 
elliptic surfaces, Kodaira discovered a number
of configurations of rational curves related to what is now
called an affine Dynkin diagram. 
The nodes of these diagrams represent rational
curves, with multiplicities marked incide the circle, and the edges
are the intersection points. We reproduce
the diagrams which are relevant to our
case. This is $\tilde E_6$, the affine Dynkin $E_6$:

\hfill

\begin{tikzpicture}[
  scale=0.7,
  every node/.style={draw, circle, minimum size=4mm, fill=cyan!20}
]
  % Center (shifted left)
  \node (c) at (-0.5,0) {3};

  % Left arm
  \node (l1) at (-2,0) {2};
  \node (l2) at (-3.5,0) {1};

  % Right arm
  \node (r1) at (1,0) {2};
  \node (r2) at (2.5,0) {1};

  % Vertical arm
  \node (u1) at (-0.5,1.5) {2};
  \node (u2) at (-0.5,3) {1};

  % Edges
  \draw (c)--(l1)--(l2);
  \draw (c)--(r1)--(r2);
  \draw (c)--(u1)--(u2);
\end{tikzpicture}

\hfill

This is  $\tilde E_7$, the affine Dynkin $E_7$:

\begin{tikzpicture}[
  scale=0.7,
  every node/.style={draw, circle, minimum size=4mm, fill=cyan!20}
]
  % Horizontal chain
  \node (n1) at (-4.5,0) {1};
  \node (n2) at (-3,0) {2};
  \node (n3) at (-1.5,0) {3};
  \node (n4) at (0,0) {4};
  \node (n5) at (1.5,0) {3};
  \node (n6) at (3,0) {2};
  \node (n7) at (4.5,0) {1};

  % Vertical branch (at the 3rd node)
  \node (v) at (0,1.5) {2};

  % Edges
  \draw (n1)--(n2)--(n3)--(n4)--(n5)--(n6)--(n7);
  \draw (n4)--(v);
\end{tikzpicture}

\hfill

This is  $\tilde E_8$, the affine Dynkin $E_8$:

\begin{tikzpicture}[
  scale=0.8,
  every node/.style={draw, circle, minimum size=4mm, fill=cyan!20}
]
  % Horizontal chain
  \node (n1) at (-6,0) {2};
  \node (n2) at (-4.5,0) {4};
  \node (n3) at (-3,0) {6};
  \node (n4) at (-1.5,0) {5};
  \node (n5) at (0,0) {4};
  \node (n6) at (1.5,0) {3};
  \node (n7) at (3,0) {2};
  \node (n8) at (4.5,0) {1};

  % Vertical branch (correct position)
  \node (v) at (-3,1.5) {3};

  % Edges
  \draw (n1)--(n2)--(n3)--(n4)--(n5)--(n6)--(n7)--(n8);
  \draw (n3)--(v);
\end{tikzpicture}

\hfill

The affine Dynkin diagram $\tilde A_{n}$
(with $n=7$ nodes on the picture) is a circle
formed by rational curves of multiplicity 1,
which intersect sequentially:

\begin{tikzpicture}[
  scale=0.7,
  every node/.style={draw, circle, minimum size=5mm, fill=cyan!20}
]
  \def\n{7} % number of nodes = 2n+1 (must be odd)

  % Nodes on a circle
  \foreach \i in {1,...,\n} {
    \node (N\i) at ({360/\n * (\i-1)}:2.5) {1};
  }

  % Edges (cycle)
  \foreach \i in {1,...,\n} {
    \pgfmathtruncatemacro{\j}{mod(\i,\n)+1}
    \draw (N\i)--(N\j);
  }
\end{tikzpicture}

\hfill

The affine  Dynkin diagram $\tilde D_n$, with $n+1$ nodes:

\hfill

\begin{tikzpicture}[
  scale=0.8,
  every node/.style={draw, circle, minimum size=4mm, fill=cyan!20},
  dot/.style={draw=none, fill=none}
]
  % Left fork (now symmetric to right)
  \node (lu) at (0,1.2) {1};
  \node (ld) at (0,-1.2) {1};
  \node (l) at (1.5,0) {2};

  % Chain with ellipsis
  \node (a1) at (3,0) {2};
  \node[dot] (dots) at (4.5,0) {$\cdots$};
  \node (a2) at (6,0) {2};

  % Right fork
  \node (r) at (7.5,0) {2};
  \node (ru) at (9,1.2) {1};
  \node (rd) at (9,-1.2) {1};

  % Edges
  \draw (l)--(lu);
  \draw (l)--(ld);
  \draw (l)--(a1);
  \draw (a1)--(dots);
  \draw (dots)--(a2);
  \draw (a2)--(r);
  \draw (r)--(ru);
  \draw (r)--(rd);
\end{tikzpicture}

\hfill

%%%%%%%%%%%%%%%%%%%%%%%%%%%%%%%%%%%%%%%%%%%%%%%%%%%%%%%%%%%%
\theorem\label{_HO_classification_Theorem_}
Let $M$ a holomorphically symplectic
manifold and $\pi:\; M \to B$ a proper Lagrangian
fibration with fibers of Fujiki class C. 
Denote by $D$ the discriminant (set of 
critical values) of $\pi$;\footnote{By 
\cite[Proposition 3.1]{_Hwang_Oguiso:characteristic_},
$D$ is a divisor in $B$.} Assume that 
$\pi$ has multiplicity 1 in codimension 1,
and let $z\in D$ be a general point. 
Denote by $\Theta$ a charateristic cycle
in $\pi^{-1}(z)$. Then $\Theta$ is
one of the singular fibers in a relative
minimal elliptic fibration listed by Kodaira
in \cite[Theorem 6.2]{_Kodaira:elliptic_} or
$A_\infty$, which is an infinite chain of
rational curves, intersecting sequentially,
and the intersection is transversal, that is,
of multiplicity 1. 

\proof \cite[Proposition 2.3]{_Hwang_Oguiso:multiple_}.
\endproof

\hfill

\remark The singular fibers in a 
minimal elliptic fibration (\cite[Theorem 6.2]{_Kodaira:elliptic_})
can be classified in terms of the affine Dynkin diagrams 
$\hat A_n$, $\hat D_n$, $\hat E_6$, $\hat E_7$, $\hat E_8$.
Each node in these diagrams corresponds to 
a characteristic leaf, with multiplicity which
is determined by the diagram, and the
edges correspond to intersections of the
characteristic lines, which are transversal.

\hfill

For further use, we reproduce Kodaira's list from
\cite[Theorem 6.2]{_Kodaira:elliptic_},
see also \cite[\S V.7]{_BHPV_}.

\hfill

\proposition
Kodaira's list of non-multiple fibers
in \cite[Theorem 6.2]{_Kodaira:elliptic_}:
\begin{description}
\item[I${}_n$] 
The affine Dynkin diagram $\tilde A_{n}$, 
pictured above. The case I${}_1$
corresponds to a rational curve with one
ordinary double point.
\item[II] A rational curve with 
one simple cusp.
\item[III] Two rational curves which
intersect in one point of even multiplicity.
\item[IV]  3 rational lines intersecting in one point
\item[I${}^*_m$, $m\geq 1$]  affine Dynkin diagram $\tilde D_{m+4}$.
\item[II${}^*$]  affine Dynkin diagram $\tilde E_{8}$.
\item[III${}^*$]  affine Dynkin diagram $\tilde E_{7}$.
\item[IV${}^*$]  affine Dynkin diagram $\tilde E_{6}$.
\end{description}
\endproof

\hfill

%%%%%%%%%%%%%%%%%%%%%%%%%%%%%%%%%%%%%%%%%%%%%%%%
\remark\label{_torus_HO_Remark_}
Let $\pi:\; M \to X$ be a Lagrangian fibration with general
fiber an abelian variety (in particular, connected), and
discriminant $D$.
By \cite[Proposition 2.3]{_Hwang_Oguiso:multiple_},
there is a torus action on a general fiber 
of $\pi:\; \pi^{-1}(D) \to D$, and this
action is transitive on characteristic curves.
This implies, in particular, that the
set of irreducible components of $\pi^{-1}(z)$,
for $z\in D$ a generic point, is a quotient
of the intersection graph of characteristic
leaves inside a characteristic curve.
Notably, in the case $A_\infty$ (which can
be realized, as shown in 
\cite[Proposition 4.13]{_Hwang_Oguiso:characteristic_}),
the set of irreducible components is arranged in a circle
which is obtained as a quotient of the
graph $A_\infty$ by $\Z$ acting by translations.

\hfill

%%%%%%%%%%%%%%%%%%%%%%%%%%%%%%%%%%%%%%%%%%%%%%%%
\claim \label{_irred_comp_HO_Claim_}
(see also 
\cite[Proposition 3.3, Proposition 3.4]{_Hwang_Oguiso:multiple_})\\
 Let $\pi:\; M \to X$ be a Lagrangian fibration with general
fiber an abelian variety (in particular, connected), and
$z\in D$ a general point, corresponding to a multiple fiber. 
Let $R$ be the incidence 
graph of irreducible components of $\pi^{-1}(z)$.
Then $R$ coincides with the graph of $R_0$ characteristic
cycles described in \ref{_HO_classification_Theorem_},
unless $R_0$ is a cycle of rational curves (called
$I_{2m}$ in Kodaira classification, and 
corresponding to the affine Dynkin diagram $\tilde A_{2m}$)
or $A_\infty$. In these two cases, $R$ is a cycle,
obtained as a quotient of $R_0$ by translations.

\hfill

\proof
Consider the map $\phi:\; R \to R_0$ taking a 
characteristic curve to its irreducible component.
By \cite[Proposition 2.3]{_Hwang_Oguiso:multiple_},
a torus $T$ acts\footnote{We repeat its definition 
later in \ref{_second_cohomo_fiber_Proposition_}, Step 2.}
on the set of characteristic cycles
transitively. This implies that $\phi$ is surjective.
By \ref{_HO_classification_Theorem_},
the characteristic cycle consists of rational
curves, unless it is just one elliptic curve.
Therefore, we can assume that characteristic
cycle is a union of rational curves.
Since the stabilizer $\St_\Theta$ of a characteristic cycle
$\Theta$ is a subgroup in $T$, it is a finite cyclic group,
acting on $\pi^{-1}(z)$ without fixed points. 
Any automorphism of a rational curve preserves
a point, hence $\St_\Theta$ acts on $R$ without
fixed points. This excludes all graphs obtained
in \ref{_HO_classification_Theorem_}, except 
$\hat A_{2n}$, $A_\infty$ and $\hat D_{2m +1}$,
because any automorphism of the rest has a fixed vertex.
The graph $\hat D_{2m +1}$ has a number of automorphisms
without fixed vertices, but they all preserve the middle
edge, which corresponds to a nodal point, hence
the corresponding action on $\Theta$ also has a fixed point.
\endproof

%%%%%%%%%%%%%%%%%%%%%%%%%%%%%%%%%%%%%%%%%%%%%%%%%%%%%%%%%%%%
\subsection{Second cohomology of the general singular fiber} 
\label{_H^2_Subsection_}
%%%%%%%%%%%%%%%%%%%%%%%%%%%%%%%%%%%%%%%%%%%%%%%%%%%%%%%%%%%%

%The following question is often referred as ``Witten's conjecture''
%(\cite???). Let $\pi:\; X \to Y$ be the Lagrangian fibration
%associated with a Hitchin system, $\dim_\C X=2n$. Is the 
%natural map $H^{2n}_c(X) \to H^{2n}(X)$ an isomorphism?
%Originally it was asked about the pairing in $L^2$-cohomology
%of $X$. Since $H^{2n}_c(X)$ is dual to $H^{2n}(X)$
%which is dual to $H_{2n}(X)$, this can be interpreted
%as a question about the intersection pairing in 
%middle homology. In the present section, we 
%compute this pairing for a neighbourhood
%of a fiber over a general point in the
%discriminant, and apply this result
%to compute the homology of the preimage 
%of a 1-dimensional disk in a base of a Lagrangian fibration.
%
%\hfill

%%%%%%%%%%%%%%%%%%%%%%%%%%%%%%%%%%%%%%%%%%%%%%%%%%%%%%%%%%%%
\proposition\label{_second_cohomo_fiber_Proposition_}
Let $\pi:\; M \to \C P^n$ be a Lagrangian fibration
on a hyperk\"ahler manifold, and $z\in D$ 
a general point.  Let $W_1,..., W_m$ be the irreducible
components of $W:=\pi^{-1}(z)$. 
Assume that the characteristic leaves of $W_i$ are
rational curves.  Denote by $S$ the incidence graph of
$W_i$ (by \ref{_irred_comp_HO_Claim_}, 
it is a tree for all affine Dynkin diagrams
except $\hat A_n$ and $A_\infty$, and homotopy equivalent
to a circle otherwise).
\begin{description}
\item[(i)] For each $W_i$, we have
$H^2(W_i, \R)=H^2(\Alb(W_i), \R)\oplus \langle
  l_i\rangle$,
where $l_i$ is a class in $H^2(W_i)$ which restricts a 
generator of the characteristic leaf $L_i$, and
$\Alb(W_i)$ the Albanese variety.
\item[(ii)] Each connected component of the 
intersection $W_i\cap W_{i+1}$ is a torus, 
and the Albanese map $\Alb:\; W_i\cap W_{i+1}\to \Alb(W_i)$
is an isogeny. This identifies the
the subspace $\Alb^*(H^*(\Alb(W_i))) \subset H^*(W_i)$
for all $W_i$.\footnote{By abuse of notation, we use
$\Alb$ for the Albanese variety and for the Albanese map.} 
Let $V:= \Alb(W_1)$. We will consider $H^2(V)$
as a subspace in $H^2(W_i)$ from now on.

\item[(iii)] 
\begin{equation}\label{_second_cohomo_when_tree_Equation_}
H^2(W)= H^2(V)\oplus \bigoplus_{i=1}^m  \langle  l_i\rangle
\end{equation} when
$S$ is a tree. When $S$ is  homotopy equivalent
to a circle, we have 
\begin{equation}\label{_second_cohomo_when_circle_Equation_}
 H^2(W)= H^2(V)\oplus 
  H^1(V)\oplus\bigoplus_{i=1}^m \langle  l_i\rangle.
 \end{equation}
\end{description}
\pstep
By \cite[Proposition 2.2]{_Hwang_Oguiso:multiple_},
the Albanese map $W_i\to \Alb(W_i)$ is a $\C P^1$-fibration.
Then the Leray spectral sequence implies
$H^2(W_i, \R)=H^2(\Alb(W_i), \R)\oplus \langle l_i\rangle$ 
(the differentials in the Leray spectral sequence
vanish for any K\"ahler fibration, as follows from
\cite{_Deligne:degenerate_}). This proves
\ref{_second_cohomo_fiber_Proposition_} (i).

\hfill

{\bf Step 2:} 
Let $T\subset \Aut(W)$ be the group generated by
the Hamiltonian vector fields associated with
the pullbacks of the holomorphic functions on $\C P^n$
defined in a neighbourhood of $z$.
By \cite[Proposition 2.3]{_Hwang_Oguiso:multiple_},
$T$ acts on transitively on the set of characteristic
curves. Since the characteristic leaves in
adjancent $W_i$, $W_{i+1}$ intersect transversally,
$T$ acts transitively on each irreducible
component of $W_i\cap W_{i+1}$. This implies
that this intersection is a torus. Applying
\cite[Proposition 2.2]{_Hwang_Oguiso:multiple_} again,
we obtain that the Albanese map $W_i \to \Alb(W_i)$
restricted to $W_i\cap W_{i+1}$ is an isogeny.
This takes care of \ref{_second_cohomo_fiber_Proposition_} (ii).

\hfill

{\bf Step 3:} In Step 3, we treat the case when 
$S$ is a tree. Then the restriction map
$H^2(W_i) \to H^2(W_i \cap W_{i+1})$ is surjective
with kernel $\langle l_i\rangle$.
Using induction on the number of vertices
in a connected subtree $S_1 \subset S$, we can assume that 
the isomorphism \eqref{_second_cohomo_when_tree_Equation_}
holds for the union $W'$ of $W_m$ enumerated by the
vertices of $S_1$, and $W= W' \cup W_m$. The Mayer-Vietoris exact
sequence gives
\begin{multline}\label{_MV_long_Equation_}
H^1(W') \oplus H^1(W_m) \to
 H^1(W_m \cap W')\stackrel a 
\to \\ \stackrel a \to H^2(W) \to H^2(W') \oplus H^2(W_m) \to H^2(W_m \cap W').
\end{multline}
Each connected component of the 
intersection $W_m \cap W'$ is a torus (Step 2)
and the map $H^1(W_m)\to  H^1(W_m \cap W')$
is surjective by the same argument with degeneration
of the Leray spectral sequence as used in Step 1.
Therefore, $a=0$  and
\eqref{_MV_long_Equation_} is reduced to
\begin{equation}\label{_MV_short_Equation_}
 0 \to H^2(W) \to H^2(W') \oplus H^2(W_m) \to H^2(W_m \cap W').
\end{equation}
The induction step gives 
$H^2(W')= H^2(V)\oplus \bigoplus_{i=1}^{m-1} \langle  l_i\rangle$,
Step 1 gives $H^2(W_m)=H^2(V)\oplus \langle  l_n\rangle$,
and Step 2 gives $H^2(W_m \cap W')=H^2(V)$.
Therefore, the exact sequence 
\eqref{_MV_short_Equation_}
can be rewritten as
\[
0 \to H^2(W) \to H^2(V)\oplus 
H^2(V)\oplus\bigoplus_{i=1}^{m} \langle  l_i\rangle \to 
H^2(V),
\]
and this immediately implies \eqref{_second_cohomo_when_tree_Equation_}.

\hfill

{\bf Step 4:} It remains to prove 
\eqref{_second_cohomo_when_circle_Equation_}.
Let $W_1, ..., W_m$ be the irreducible components
of $W$ arranged in a circle, and $W'= \bigcup_{i=1}^{m-1} W_i$.
By Step 3, 
$H^2(W')= H^2(V)\oplus \bigoplus_{i=1}^{m-1} \langle  l_i\rangle$.
Consider the Mayer-Vietoris exact sequence for
$W=W'\cup W_m$,
\begin{multline}\label{_MV_for_circle_long_Equation_}
H^1(W)\to H^1(W') \oplus H^1(W_m) \to
 H^1(W_m \cap W') \to \\
\to H^2(W) \to H^2(W') \oplus H^2(W_m) \to H^2(W_m \cap W').
\end{multline}
However, unlike in Step 3, the intersection
$W_m \cap W'$ has two connected components, which
are tori, isogeneous to $\Alb(W_i)$.
The restriction map $H^1(W)\to H^1(W') \oplus H^1(W_m)$
takes $H^1(V)$ to $H^1(V)\oplus H^1(V)$,
hence its cokernel is $H^1(V)$, which is injectively
mapped to $H^1(W_m \cap W') =H^1(V)\oplus H^1(V)$.
This reduces the long exact sequence \eqref{_MV_for_circle_long_Equation_}
to a shorter one
\begin{equation}\label{_MV_for_circle_short_Equation_}
0\to H^1(V)
\to H^2(W) \to H^2(W') \oplus H^2(W_m) \to H^2(W_m \cap W').
\end{equation}
Now, the same argument as in Step 3 is used to transform
\eqref{_MV_for_circle_short_Equation_} to 
\eqref{_second_cohomo_when_circle_Equation_}.
\endproof

\hfill

{\bf Acknowledgments.} The second named author thanks IMPA for their 
hospitality in January 2020 when part of this project began. 
We are grateful to Claire Voisin and Andrey Soldatenkov who
helped us to find the error in the first version of this paper.
Many thanks to Jason Starr for references to rational 
connectedness in relation to rational singularities.
We are grateful to J. Koll\'ar and G. Sacc\`a
for their insightful comments and email exchange,
and to the anonymous referee for his or her attention
and suggestions.

\small
\noindent {\sc Fedor A. Bogomolov\\
Department of Mathematics\\
Courant Institute, NYU \\
251 Mercer Street \\
New York, NY 10012, USA,} \\
\tt bogomolov@cims.nyu.edu,\\
{\sc also: \\
 Laboratory of Algebraic Geometry,\\
National Research University Higher School of Economics,\\
Department of Mathematics, 6 Usacheva street, Moscow, Russia.}
\\

\noindent {\sc Ljudmila Kamenova\\
Department of Mathematics, 3-115 \\
Stony Brook University \\
Stony Brook, NY 11794-3651, USA,} \\
\tt kamenova@math.sunysb.edu
\\

\noindent {\sc Misha Verbitsky\\
            {\sc Instituto Nacional de Matem\'atica Pura e
              Aplicada (IMPA) \\ Estrada Dona Castorina, 110\\
Jardim Bot\^anico, CEP 22460-320\\
Rio de Janeiro, RJ - Brasil}\\
\tt verbit@impa.br
}

\end{document}